\title{ ~~\\ On the distribution of the order and index of
$g({\rm mod~}p)$ over residue classes}
\author{Pieter Moree}
\documentclass[12pt]{article}
\usepackage{amssymb, latexsym, amsfonts}
\def\@ptsize{2}
\newtheorem{Thm}{Theorem}

\newtheorem{Lem}{Lemma}
\newtheorem{Cor}{Corollary}

\newtheorem{Prop}{Proposition}
\newcommand{\qed}{\hfill $\Box$}

\begin{document}
\date{}
\maketitle
{\def\thefootnote{}
\footnote{\noindent Max-Planck-Institute f\"ur Mathematik,
Vivatsgasse 7, D-53111 Bonn, Deutschland, E-mail:
moree@mpim-bonn.mpg.de}}
{\def\thefootnote{}
\footnote{{\it Mathematics Subject Classification (2000)}.
11N37, 11N69, 11R45}}

\begin{abstract}
\noindent For a fixed rational number $g\not\in \{-1,0,1\}$ and
integers $a$ and $d$ we consider the
set  $N_g(a,d)$ of primes $p$ for which the order of $g({\rm mod~}p)$ is
congruent to $a({\rm mod~}d)$.
For $d=4$ and $d=3$ we show that, under the
Generalized Riemann Hypothesis (GRH), these sets have a natural
density $\delta_g(a,d)$
and compute it. 
The results for $d=4$ generalise earlier work by Chinen and Murata.
The case $d=3$ was apparently not considered before.
\end{abstract}

\section{Introduction}
\label{introdu}
Let $g\not\in \{-1,0,1\}$ be a rational number (this
assumption on $g$ will be maintained throughout this paper). For $u$
a rational number, let $\nu_p(u)$ denote the exponent of $p$ in the canonical
factorisation of $u$ (throughout the letter $p$ will be used to indicate
prime numbers). If $\nu_p(g)=0$, then there exists a smallest positive integer
$k$ such that $g^k\equiv 1({\rm mod~}p)$.
We put ord$_g(p)=k$.
This number is the {\it (residual)
order} of $g({\rm mod~}p)$.
The index of the subgroup generated by $g$ mod $p$ inside the multiplicative
group of residues mod $p$, $[({\mathbb Z}/p{\mathbb Z})^\times:
\langle g({\rm mod~}p)\rangle]$, is denoted by $r_g(p)$ and
called the {\it (residual) index} mod $p$ of $g$.
Although ${\rm ord}_g(p)$ and $r_g(p)$ satisfy the easy relation
\begin{equation}
\label{indord}
{\rm ord}_g(p)r_g(p)=p-1,
\end{equation}
the functions themselves fluctuate quite irregularly. Given this it comes
perhaps not as a surprise that a simple question such as Artin's primitive
root conjecture (1927), which asserts that $\{p:r_g(p)=1\}$ is infinite if
$g$ is not a square, remains unsolved. On the assumption of the Riemann
Hypothesis for a certain class of Dedekind zeta functions, however, this
was proved by C. Hooley \cite{Hooley}. Many variations of Artin's conjecture
have been considered in the course of time, the most far reaching in
\cite{Lenstra}.\\
\indent Many authors studied the divisibility of the order by some prescribed
integer $d$. The case $d=2$ for example is closely related to the
non-divisiblity of certain integer sequences by a prescribed prime.
We say that an integer sequence $S=\{s_j\}_{j=1}^{\infty}$ is divisible by an
integer $m$, if there exists an integer $k$ such that $m|s_{k}$. It is
easy to see that for a prime $p$ with $\nu_p(g)=0$ the sequence
$S(g)=\{g^j+1\}_{j=1}^{\infty}$ is divisible by $p$ if and only if
ord$_g(p)$ is even. Hasse \cite{Hasse1, Hasse2} showed that the set of
prime divisors of the aforementioned sequence has a Dirichlet density.
It is not difficult to extend his argument to show that these sequences
have a natural density (hereafter we merely write density instead
of natural density) of prime divisors. For $g=10$ a prime $p\ne 2,5$
divides $S(g)$ if and only if the period of the decimal expansion of
$1/p$ is even, cf. \cite{Odoni}.
Using some algebraic number
theory these results can be extended to some other 
well-known sequences,
cf. \cite{Ballot, Lagarias, MoreeS}. In
all of these cases the density can be computed unconditionally
and turns out to be a rational number. For example, the density of prime
divisors of $S(2)$ is 17/24.\\
\indent Now let $d>2$ be given. By similar methods the divisibility of
the order by $d$ or the coprimality of the order with $d$ can be studied.
In this direction we especially like to  mention K. Wiertelak,
who wrote many papers on this subject, starting in the seventies of the
previous century. See \cite{Wiertelak2} for his most recent paper.
Again one can prove that the density of 
the set of such primes exists and is rational.\\
\indent In the light of the extensive literature on the case where the
order is divisible by $d$, it is somewhat surprising that the
question of how the order is distributed over the various residue classes
mod $d$ has up to this century only been considered for $d=2$. The purpose of this
paper and its sequel(s) is to address this question for various other values
of $d$. For the understanding of the general case it is in my viewpoint crucial to
first study a particular case in detail, for which we take $d=4$.\\
\indent For $d=4$ our main interest is in the set $N_g(a,4)$, but it
turns out to be fruitful to consider $N_g(1,2^s;j,4)(x)$ and
$N_g(3,4;j,4)(x)$ separately, where 
$N_g(a_1,d_1;a_2,d_2)(x)$ counts the number of  primes $p\le x$ satisfying $\nu_p(g)=0$ 
for which
$p\equiv a_1({\rm mod~}d_1)$ and ord$_g(p)\equiv a_2({\rm mod~}d_2)$.
For convenience we denote $N_g(0,1;a,d)(x)$ by $N_g(a,d)(x)$.
Although the functions
 $N_g(1,2^s;1,4)(x)$ and $N_g(1,2^s;3,4)(x)$ 
 are more complicated (see Theorem \ref{triplesum}) to describe, they 
 turn out to be asymptotically
equal under GRH. For the  more
easily describable functions $N_g(3,4;1,4)(x)$ and $N_g(3,4;3,4)(x)$ (vide Lemma 
\ref{simpeltoch}),
the asymptotic behaviour can be different.\\
\indent For $s|r$ the number field
$\mathbb Q(\zeta_r,g^{1/s})$ will be denoted by $K_{r,s}$. By
$\pi_L(x)$ we denote the number of rational primes $p\le x$ that are
unramified in the number field $L$ and split completely in $L$. As
usual we let Li$(x)$ denote the logarithmic integral, that
is Li(x)=$\int_2^x {dt\over \log t}$.\\
\vfil\eject
\indent For the primes $p\equiv 1({\rm mod~}2^s)$, $s\ge 2$, we find:
\begin{Thm}
\label{main1} Write $g=g_1/g_2$ with $g_1,g_2$ integers.
Let $s\ge 2$. For $j=0$ and $j=2$ we have
\begin{equation}
\label{densprecise}
N_g(1,2^s;j,4)(x)=\delta_g(1,2^s;j,4){\rm Li}(x)+
O\left({x(\log \log x)^4\over \log^3 x}\right),
\end{equation}
where
\begin{equation}
\label{dens0mod4}
\delta_g(1,2^s;0,4)=2^{1-s}-\sum_{r\ge s}\Big({1\over [K_{2^r,2^{r-1}}:\mathbb
Q]}-{1\over [K_{2^{r+1},2^{r-1}}:\mathbb Q]}\Big){\rm ~and~}
\end{equation}
$$\delta_g(1,2^s;2,4)=\sum_{r\ge s}\Big({1\over [K_{2^r,2^{r-1}}:\mathbb
Q]}-{1\over [K_{2^{r+1},2^{r-1}}:\mathbb Q]}-{1\over [K_{2^r,2^{r}}:\mathbb
Q]}+{1\over [K_{2^{r+1},2^{r}}:\mathbb
Q]}\Big).$$
For $j=1$ and $j=3$ we have, under GRH, 
$$N_g(1,2^s;j,4)(x)={\delta_g(1,2^s;1,2)\over 2}{x\over \log x}
+O\left(\log|g_1g_2|{x\over \log^{3/2}x}\right),$$
where
$$\delta_g(1,2^s;1,2)=\sum_{r\ge s}\Big({1\over [K_{2^r,2^{r}}:\mathbb
Q]}-{1\over [K_{2^{r+1},2^{r}}:\mathbb Q]}\Big)$$is the density of the set
$N_g(1,2^s;1,2)$ and the implied constant is absolute.
\end{Thm}
\indent For the primes $p\equiv 3({\rm mod~}4)$ we find:
\begin{Thm}
\label{main3}
Write $g=g_1/g_2$ with $g_1$ and $g_2$ integers.
Let $\psi_0,\psi_1$ denote the principal, respectively non-principal
character mod $4$. Let $h_{\psi_1}(v)=\sum_{d|v}\mu(d)\psi_1(v/d)$, where
$\mu$ denotes the M\"obius function.
We have 
$$\cases{N_g(3,4;0,4)(x)=0; \cr
N_g(3,4;2,4)(x)=\#\{p\le x:p\equiv 3({\rm mod~}4),~({g\over p})=-1\}.}$$
Assuming GRH we have, when $j$ is odd,
$$N_g(3,4;j,4)(x)={1\over 2}\#\{p\le x:p\equiv 3({\rm mod~}4),~({g\over
p})=1\}$$$$
+(-1)^{j-1\over 2}{\Delta_g\over 4}{x\over \log x}+
O\left(\log|g_1g_2|{x\over \log^{3/2}x}\right),$$
where
$$\Delta_g=\sum_{\sqrt{-2}\in K_{2v,2v}\atop 2\nmid v}
{h_{\psi_1}(v)\over [K_{2v,2v}:\mathbb Q]}
-\sum_{\sqrt{2}\in K_{2v,2v}\atop 2\nmid v}
{h_{\psi_1}(v)\over [K_{2v,2v}:\mathbb Q]}$$
and the implied constant is absolute.
\end{Thm}
It is clear that $N_g(j,4)(x)=N_g(1,4;j,4)(x)+N_g(3,4;j,4)(x)$ and
we leave it to the reader to add the estimates for the latter two
quantities given in Theorem \ref{main1}, respectively Theorem \ref{main3}
(Corollary \ref{samen} provides an example).\\
\indent In Section \ref{Explicit} we derive explicit versions of
Theorem \ref{main1} and Theorem \ref{main3}. 
For $s\ge 0$ it follows that $\delta_g(1,2^s;j,4)$ exists
and is in $\mathbb Q+\mathbb QA_{\psi_1}$, where
$$A_{\psi_1}=\prod_{p\equiv 3({\rm mod~}4)}\left(1-{2p\over
(p^2+1)(p-1)}\right)=0.643650679662525\dots$$
As an example we mention the following corollary to Theorem
\ref{main1} and Theorem \ref{main3} (for the notation $h$ and $D$
we refer to Lemma \ref{degree}; for a non-zero real number $r$ we denote its sign by sgn$(r)$).
\begin{Cor} {\rm (GRH)}.
\label{samen}
Suppose that $h=1$ and $j$ is odd. Then
$\delta_g(1,4)=\delta_g(3,4)=1/6$ unless $D$ is divisible by 8 and has
no prime divisor congruent to
$1({\rm mod~}4)$, in which case we have
$$\delta_g(j,4)=
\cases{{1\over 6}+{\rm sgn}(g)A_{\psi_1}{(-1)^{j+1\over 2}\over
8}\prod_{p|{D\over 8}}{2p\over p^3-p^2-p-1} &if $D\ne 8$;\cr
{7\over 48}+{\rm sgn}(g)A_{\psi_1}{(-1)^{j+1\over 2}\over
8} &if $D=8$.}$$
\end{Cor}
\indent In case $d=3$ similar results to those for $d=4$ are obtained
in Section \ref{modulusdrie}.
In particular we will show that, under GRH,
$\delta_g(1,3^s;j,3)$ exists for $s\ge 0$ and that
$\delta_g(1,3^s;j,3)\in \mathbb Q+\mathbb Q A_{\xi_1}$, where
$$A_{\xi_1}=\prod_{p\equiv 2({\rm mod~}3)}\left(1-{2p\over
(p^2+1)(p-1)}\right)=0.173977122429634\dots,$$
and $\xi_1$ denotes the non-principal character mod 3.
The rational numbers involved we explicitly compute.\\
\indent The analogous problem of studying the primes for which the {\tt index}
is congruent to $a({\rm mod~}d)$ turns out to be far easier 
(see Section \ref{indexmodd}). Nevertheless, at least for $d=3$ and
$d=4$, we again find, under GRH, that these densities exist and
are in $\mathbb Q+\mathbb Q A_{\xi_1}$, respectively
$\mathbb Q + \mathbb Q A_{\psi_1}$.\\
\indent Instead of requiring GRH it is enough to require that RH holds
for every field $\mathbb Q(\zeta_{r},g^{1/s})$ with $s|r$. Indeed, if
a given result is under GRH we mean that we require RH to hold for
every field that occurs in the proof of this result.\\
\indent In a sequel to this paper by a slightly different method the
case where $d$ is an odd prime power is investigated (but less
explicitly), see \cite{MorII}.\\
\indent The density $\delta_g(j,4)$ with $g$ a positive integer
that is not a pure power (i.e. $h=1$ in the notation of
Lemma \ref{degree}), was first studied by Chinen and Murata
in \cite{CM1,CM2,CM3,CMnogeen, CM4,CM5}, culminating (in \cite{CM5}) in
their proof of Corollary 1 for the case $g>0$.

\section{Preliminaries}
\subsection{The index and algebraic number theory}
In this section we recall some well-known arguments from the theory of primitive roots
that are essential for an understanding of the rest of this paper.\\
\indent The index can be easily related to algebraic number theory and
 by using (\ref{indord}) we then can get a grip on the order.
Thus a prime $p$ that satisfies $k|r_g(p)$ must obviously satisfy
$p\equiv 1({\rm mod~}k)$ and $g^{(p-1)/k}\equiv 1({\rm mod~}p)$, in other
words it must split completely in the field $\mathbb Q(\zeta_{k},
g^{1/k})(=K_{k,k})$. On the other hand a prime $p$ that satisfies the latter
condition satisfies $k|r_g(p)$. Then, by the principle of inclusion
and exclusion, we can describe for example the set of primes $p$ that
satisfy $r_g(p)=k$. Note that $r_g(p)=k$ iff
$k|r_g(p)$ and $qk\nmid r_g(p)$ for any prime $q$.
Let $R_g(a,f;t)$ denote the set of primes $p$ with
$p\equiv a({\rm mod~}f)$ and $r_g(p)=t$. Let $R_g(a,f;t)(x)$ denote
the number of primes $p\le x$ in $R_g(a,f;t)$.
Using the
principle of inclusion and exclusion, we then find that
\begin{equation}
\label{inclusionex}
R_g(a,f;t)(x)=\sum_{n=1}^{\infty}\mu(n)\#\{p\le x:~p\equiv a({\rm mod~}f),~
(p,K_{nt,nt}/\mathbb Q)=id\}.
\end{equation}
By
$(p,K/\mathbb Q)$ we denote the Frobenius symbol. 
We have $(p,K/\mathbb Q)=id$ iff $p$ is unramified and splits
completely in $K$.
Since sets of the
form $\{p:p\equiv a({\rm mod~}f),~
(p,K_{r,n}/\mathbb Q)=id\}$, will occur rather frequently in the sequel,
we will denote them by $S_g(a,f;r,n)$ and the corresponding counting
function by $S_g(a,f;r,n)(x)$.
Assuming GRH, it follows from
\cite{Lenstra} that $R_g(a,f;t)$ has a density.\\
\indent Sofar this
density has only been evaluated in terms of
an Euler product (singular series) in the case $t=1$ with $a$ and
$f$ arbitrary \cite{LMS}, or in the case $f|2$ and $t$ arbitrary
\cite{Murata, Wagstaff}. For example, for $t=1$ and $2|f$ the density,
under GRH, is a rational multiple of the Artin constant
$\prod_p(1-{1\over p(p-1)})$ \cite{Hooley}.
For an unified Galois theoretic treatment of finding Euler products for
these cases see \cite{LMS}.\\
\indent For our purposes such
an evaluation of the density of $R_g(a,f;t)$ will, however,
be irrelevant, an evaluation
in terms of an infinite series
will be sufficient. The tool to arrive
at such an expression for the density is the Chebotarev
density theorem:
\begin{Thm}
\label{chebotarev}
{\rm (GRH)}. Let $K$ be an algebraic number field, let $L/K$ be a finite
Galois extension and $C$ be a conjugacy class in $G={\rm Gal}(L/K)$. We
let $\pi(x;L/K,C)$ denote the number of unramified prime ideals 
$\frak{p}$ in $K$
such that $(\frak{p},L/K)=C$ and $N\frak{p}\le x$. Then, under RH
for the field $L$ we have
\begin{equation}
\label{chebotarevgrh}
\pi(x;L/K,C)={\#C\over \# G}{\rm Li}(x)+O\left({\# C\over \#
G}\sqrt{x}\log(d_Lx^{[L:\mathbb Q]})\right),~{\rm as~}x\rightarrow\infty,
\end{equation}
where $d_L$ denotes the discriminant of $L$.
\end{Thm}
{\tt Remark}. The proof of Theorem \ref{chebotarev} is in essence
due to Lagarias and Odlyzko, the present
formulation is due to Serre \cite[p. 133]{Serre}, who
removed `un terme parasite' in the formulation of Lagarias
and Odlyzko (\cite[Theorem 1.1]{LO}). In case $C=id$
the result was proved earlier by Lang \cite{Lang}. For several
variants of Artin's primitive root conjecture Lang's result is
all one needs. There are also unconditional variants that certainly allow 
us to deduce
that $\pi(x;L/K,C)\sim {\#C\over \# G}{\rm Li}(x)$, as $x$ tends
to infinity.\\

\indent For an arbitrary integer
$m\ge 1$ let us see how Chebotarev's density theorem can be used to estimate
$S_g(a,f;r,n)(x)$, where $n|r$.
To this end we consider the compositum of the
fields $\mathbb Q(\zeta_f)$ and $K_{r,n}$, that is
$K_{[f,r],n}$, where by $[f,r]$ we denote the lowest common multiple
of $f$ and $r$.
Let $K_1$ and $K_2$ be number fields that are Galois.
If there is an automorphism
$\sigma_1\in {\rm Gal}(K_1/\mathbb Q)$ and an automorphism
$\sigma_2\in {\rm Gal}(K_2/\mathbb Q)$ such that $\sigma_1=\sigma_2$ on
$K_1\cap K_2$, then there is an unique $\sigma\in K_1\cdot K_2$, the
compositum of $K_1$ and $K_2$ such that $\sigma|_{K_j}=\sigma_j$ for
$j=1,2$. Now in order to apply the Chebotarev density theorem, we
have to count the number of 
elements in the conjugacy class
of $\sigma\in {\rm Gal}(K_{[f,r],n}/\mathbb
Q)$, where $\sigma$ is such that $\sigma|_{\mathbb Q(\zeta_f)}=\sigma_{a,f}$, where
$\sigma_{a,f}\in {\rm Gal}(\mathbb Q(\zeta_f)/\mathbb Q)$ is uniquely
determined by $\sigma_{a,f}(\zeta_f)=\zeta_f^a$, and
$\sigma_{a,f}|_{K_{r,n}}=id$. By the above remark such a $\sigma$ exists,
and is unique, if and only if
$\sigma|_{\mathbb Q(\zeta_f)\cap K_{r,n}}=id$. 
Note that a conjugate $\tau \sigma \tau^{-1}$ acts trivially on
$K_{r,n}$ and can be regarded as an element of Gal$(K_{[f,r],n}/K_{r,n})$, 
which is a subgroup of the abelian 
group $(\mathbb Z/f\mathbb Z)^*$. Hence $\tau \sigma
\tau^{-1}$ acts as $\sigma$. We conclude that the conjugacy class
has one element if $\sigma|_{\mathbb Q(\zeta_f)\cap K_{r,n}}=id$ and
zero otherwise.
By this argument we
expect from (\ref{inclusionex}), assuming there is enough cancellation in
the error terms, that the density of $R_g(a,f;t)$ is given by
\begin{equation}
\label{lenstraoud}
\lim_{x\rightarrow \infty}{R_g(a,f;t)(x)\over
\pi(x)}=\sum_{n=1}^{\infty}{\mu(n)c_1(a,f,nt)\over [K_{[f,nt],nt}:\mathbb Q]},\end{equation}
where
$$c_1(a,f,nt)=\cases{1 &if
$\sigma_{a,f}|_{\mathbb Q(\zeta_f) \cap 
K_{nt,nt}}=id$;\cr
0 &otherwise,}$$
and $\pi(x)=\sum_{p\le x}1$.
By \cite{Lenstra} we know that (\ref{lenstraoud}) holds true, under GRH.\\

\subsection{Field degrees and intersections}
In order to explicitly evaluate certain densities in this paper, the 
following result will play a crucial role. The notations $D$, $g_0$,  
$h$ and $n_r$ will reappear again and again in the sequel. If
$a$ and $b$ are integers, then by $(a,b)$ and $[a,b]$ we denote
the greatest common divisor, respectively lowest common multiple of
$a$ and $b$.
\begin{Lem}
\label{degree}
Write $g=\pm g_0^h$, where $g_0$ is
positive and not an
exact power of a rational. 
Let $D$ denote the discriminant of the field $\mathbb Q(\sqrt{g_0})$.
Put $m=D/2$ if $\nu_2(h)=0$ and $D\equiv 4({\rm mod~}8)$
or $\nu_2(h)=1$ and $D\equiv 0({\rm mod~}8)$, and
$m=[2^{\nu_2(h)+2},D]$ otherwise. 
Put $$n_r=\cases{m &if $g<0$ and $r$ is odd;\cr
[2^{\nu_2(hr)+1},D] &otherwise.}$$
We have
$$[K_{kr,k}:\mathbb Q]=[\mathbb Q(\zeta_{kr},g^{1/k}):\mathbb Q]={\varphi(kr)k\over
\epsilon(kr,k)(k,h)},$$
where, for $g>0$ or $g<0$ and $r$ even we have
$$\epsilon(kr,k)=\cases{2 &if $n_r|kr$;\cr
1 &if $n_r\nmid kr$,}$$
and for $g<0$ and $r$ odd we have
$$\epsilon(kr,k)=\cases{2 &if $n_r|kr$;\cr
 {1\over 2} &if $2|k$ and $2^{\nu_2(h)+1}\nmid k$;\cr
1 &otherwise.}$$
\end{Lem}
{\it Proof}. For $r=1$ the
result follows from Proposition 4.1 of \cite{Wagstaff} (see also the proof
of Theorem 2.2 in \cite{Wagstaff}).  For $r>1$ the result follows from the
case where $r=1$ on noting that, with ${\tilde g}=g^r$, we have
$\mathbb Q(\zeta_{kr},g^{1/k})=\mathbb Q(\zeta_{kr},{\tilde
g}^{1/kr})$. The distinction between $r$ even and $r$ odd arises in the
case $g<0$ since ${\tilde g}$ is then positive or negative, according to
whether $r$ is even or odd, respectively.\qed\\

In our analytic considerations we need an upper bound for the discriminant
of the field $K_{kr,k}$.
\begin{Lem}
\label{discriminantje}
The discriminant $D'$ of the field $K_{rk,k}$ satisfies
$$\log|D'|\le rk\Big(\log(rk)+\log(|g_1g_2|)\Big),$$
where $g_0=g_1/g_2$ and $g_1$ and $g_2$ are integers.
\end{Lem}
{\it Proof}. If $L_1/\mathbb Q$ and $L_2/\mathbb Q$ are two extension fields and
$L$ is their compositum, then the associated discriminants (over $\mathbb Q$) 
satisfy $d_{L}|d_{L_1}^{[L:L_1]}d_{L_2}^{[L:L_2]}$.
From this we have the estimate
\begin{equation}
\label{discri}
\log |d_{L}|\le [L_2:\mathbb Q]\log|d_{L_1}|+[L_1:\mathbb Q]\log|d_{L_2}|.
\end{equation}
It is well-known that the discriminant of the cyclotomic field
$\mathbb Q(\zeta_m)$ and the field $\mathbb Q(g^{1/n})$ divide
$m^{\varphi(m)}$, respectively $(ng_1g_2)^n$ (see e.g. \cite{Birch}).
On invoking these estimates, the
result then follows from (\ref{discri}). \qed\\

From cyclotomy we recall the following well-known result.
\begin{Lem}
\label{conductor}
We have $\mathbb Q(\sqrt{g})\subseteq \mathbb Q(\zeta_f)$ iff
$\Delta|f$, where $\Delta$ denotes the discriminant
of the field $\mathbb Q(\sqrt{g})$.
\end{Lem}
{\it Proof}. See e.g. \cite{Weiss}.\qed\\

In order to evaluate the densities for the modulus 4 we need the
following result, which can be easily deduced from 
Lemma \ref{degree} and the previous lemma.
\begin{Lem}
\label{intersection}
Put $L_v=\mathbb Q(\zeta_8)\cap \mathbb Q(\zeta_{2v},g^{1/2v})$.
Let $v$ be odd.\\
If $h$ is odd and $D\nmid 8v$, then $L_v=\mathbb Q$.\\ 
If $h$ is odd and $D|8v$, then
$$L_v=
\cases{
\mathbb Q(\sqrt{{\rm sgn}(g)}~)             &if $D\equiv 1({\rm mod~}4)$;\cr
\mathbb Q(\sqrt{{\rm sgn}(-g)}~)          &if $D\equiv 4({\rm mod~}8)$;\cr
\mathbb Q(\sqrt{2\cdot {\rm sgn}(g)}~)   &if $D\equiv 8({\rm mod~}32)$;\cr
\mathbb Q(\sqrt{2\cdot {\rm sgn}(-g)}~)  &if $D\equiv 24({\rm mod~}32)$.
}$$
If $h$ is even, then $L_v=\mathbb Q(\sqrt{{\rm sgn}(g)}~)$.
\end{Lem}
{\it Proof}. On using that
$$[\mathbb Q(\zeta_{8v},g^{1/2v}):\mathbb Q]=4{[\mathbb Q(\zeta_{2v},g^{1/2v}):\mathbb Q]\over [L_v:\mathbb Q]},$$
it follows by Lemma \ref{degree} that $[L_v:\mathbb Q]=\epsilon(2v,4)/\epsilon(2v,1)$. Using
Lemma \ref{degree} it follows that in each case the field claimed
to equal $L_v$ has the correct degree.\\
\indent Let us first consider the case where $g>0$. 
We may suppose that $h$ is odd and $D|8v$, since in
the remaining cases the degree of $L_v$ is 1 and hence $L_v=\mathbb Q$. It remains to show that
$i\in L_v$ if $D\equiv 4({\rm mod~}8)$, $\sqrt{2}\in L_v$ if
$D\equiv 8({\rm mod~}32)$ and $\sqrt{-2}\in L_v$ if $D\equiv 24({\rm mod~}32)$.\\
i) $D\equiv 4({\rm mod~}8)$. In this case the 
discriminant of $\mathbb Q(\sqrt{-g})$, $D/4$, divides $v$ and thus by
Lemma \ref{conductor}, $\sqrt{-g}\in \mathbb Q(\zeta_v)\subseteq K_{2v,2v}$. Since also $\sqrt{g}\in K_{2v,2v}$ it follows
that $i\in K_{2v,2v}$.\\
ii) $D\equiv 8({\rm mod~}32)$. Now $\sqrt{2g}\in \mathbb Q(\zeta_v)$
and thus $\sqrt{2}\in K_{2v,2v}$.\\
iii) $D\equiv 24({\rm mod~}32)$. Now $\sqrt{-2g}\in \mathbb Q(\zeta_v)$
and thus $\sqrt{-2}\in K_{2v,2v}$.\\
\indent Suppose $g<0$. If $h$ is even, we have to
show that $i\in \mathbb Q(\zeta_{2v},g^{1/2v})=\mathbb Q(\zeta_{2v},\zeta_{4v}g_0^{h/2v})$.
Since $(\zeta_{4v}g_0^{h/2v})^v\in \mathbb Q(i)\backslash \mathbb Q$, this is clear.
As before we may now suppose that $h$ is odd and $D|8v$.
It remains to show that
$i\in L_v$ if $D\equiv 1({\rm mod~}4)$, $\sqrt{-2}\in L_v$ if
$D\equiv 8({\rm mod~}32)$ and $\sqrt{2}\in L_v$ if 
$D\equiv 24({\rm mod~}32)$. Note that
$(\zeta_{4v}g_0^{h/2v})^v$ is a rational multiple of $\sqrt{-g_0}$ and
hence $\sqrt{-g_0}\in K_{2v,2v}$.\\
i) $D\equiv 1({\rm mod~}4)$. We have $\sqrt{g_0}\in \mathbb Q(\zeta_{v})$ and
since $\sqrt{-g_0}\in K_{2v,2v}$, it
follows that $i\in K_{2v,2v}$.\\
ii) $D\equiv 8({\rm mod~}32)$. We have $\sqrt{2g_0}\in \mathbb Q(\zeta_{v})$
and since $\sqrt{-g_0}\in K_{2v,2v}$, it follows that $\sqrt{-2}\in L_v$.\\
iii) $D\equiv 24({\rm mod~}32)$. We have 
$\sqrt{-2g_0}\in \mathbb Q(\zeta_{v})$
and since $\sqrt{-g_0}\in K_{2v,2v}$, it follows that $\sqrt{2}\in L_v$. \qed\\

\noindent Lemma \ref{intersection} 
allows one to establish the following property of $\Delta_g$.
\begin{Prop}
\label{prop1}
We have $\Delta_{-g}=-\Delta_g$.
\end{Prop}
{\it Proof}. If $h$ is even, then $\Delta_{-g}=\Delta_g=0$. If $h$ and
$v$ are odd and $8|D$ then 
$[\mathbb Q(\zeta_{2v},g^{1/2v}:\mathbb Q]=
[\mathbb Q(\zeta_{2v},(-g)^{1/2v}:\mathbb Q]$ 
by Lemma \ref{degree}. The result now
follows easily on invoking Lemma \ref{intersection}. \qed\\

\noindent 

\subsection{Index $t$ revisited}
In this section we extend some results of Murata \cite{Murata}, which he
established for squarefree integers $\ge 2$, to arbitrary
$g\in \mathbb Q\backslash\{-1,0,1\}$. However, the method of proof employed in this 
paper is rather different.\\
\indent We will need that $\sum_{r>x}1/(r\varphi(r))=O(1/x)$ (for
a proof see e.g. \cite[Lemma 8.4]{Matthews}), the following result shows
that this order of growth cannot be improved upon.
\begin{Lem} 
\label{goldbach}
We have
$$\sum_{n\le x}{1\over n\varphi(n)}
=\prod_p\left(1+{p\over (p^2-1)(p-1)}\right)-{315\zeta(3)\over 2\pi^4x}
+O\left({\log x\over x^2}\right).$$
\end{Lem}
{\it Proof}. Landau \cite{Landau} has shown that
\begin{equation}
\label{landau}
\sum_{n\le x}{1\over \varphi(n)}={315\zeta(3)\over 2\pi^4}
\{\log x+\gamma-\sum_p{\log p\over p^2-p+1}\}
+O\left({\log x\over x}\right),
\end{equation}
where $\gamma$ denotes Euler's constant.
Using (\ref{landau}) and the Euler identity, the result follows on
partial integration. Alternatively one can apply Landau's method for
establishing (\ref{landau}) to the sum  $\sum_{n\le x}1/(n\varphi(n))$. \qed\\

The estimate
\begin{equation}
\label{opzijnkop}
{1\over [K_{vt,vt}:\mathbb Q]}\le {2h\over vt\varphi(vt)},
\end{equation}
ensures that
$$A(g,t):=\sum_{v=1}^{\infty}{\mu(v)\over [K_{vt,vt}:\mathbb Q]}$$
converges absolutely.
\begin{Lem}
\label{tail}
We have
$\sum_{t\le y}A(g,t)=1+O({h\over y})$, where the implied constant is absolute.
\end{Lem}
{\it Proof}. One easily checks that $\sum_{t=1}^{\infty}A(g,t)=1$.
Using (\ref{opzijnkop}), $\varphi(vt)\ge \varphi(v)\varphi(t)$ and
Lemma \ref{goldbach}, we infer that
$$\sum_{t\ge y}A(g,t)
=O\left(\sum_{t\ge y}\sum_{v=1}^{\infty}{h\over vt\varphi(vt)}\right)
=O\left(\sum_{t\ge y}{h\over t\varphi(t)}\sum_{v=1}^{\infty}{1\over v\varphi(v)}\right)
=O({h\over y}),$$
where the implied constants are all absolute. \qed

\begin{Thm} {\rm (GRH)}.
\label{precies} Write $g=g_1/g_2$, with $g_1$ and $g_2$ integers.
For $t\le x^{1/3}$ we have
$$R_g(0,1;t)(x)=A(g,t){x\over \log x}+
O\left({x\log \log x\over \varphi(t)\log^2 x}
+{x\log|g_1g_2|\over \log^{2}x}\right),$$
where the implied constant is absolute.
\end{Thm}
{\it Proof}. Since the proof is carried out along the lines of
Hooley's proof \cite{Hooley}, we only sketch it. Let
$$M_g(x,y)=\#\{p\le x: t|r_g(p),~qt\nmid r_g(p),~q\le y\}~{\rm ~and~}$$
$$M_g(x,y,z)=\#\{p\le x: qt|r_g(p),~y\le q\le z\},$$
where $q$ denotes a prime number.
Note that
$$R_g(0,1;t)(x)=M_g(x,\tau_1)+
O(M_g(x,\tau_1,\tau_2))+O(M_g(x,\tau_2,\tau_3))
+O(M_g(x,\tau_3,{x-1\over t})).$$
We take $\tau_1={\log x/6}$, $\tau_2=\sqrt{x}\log^{-2}x$
and $\tau_3=\sqrt{x}\log x$. 
We use the starting observation that
$$R_g(0,1;t)(x)=\sum_{v=1}^{\infty}\mu(v)\pi_{K_{vt,vt}}(x).$$
Using  Lemma \ref{discriminantje} and Theorem \ref{chebotarev} we
can estimate $\pi_{K_{vt,vt}}(x)$, under GRH. Proceeding as 
Hooley did, we then obtain that both $M_g(x,\zeta_1,\zeta_2)$ and
$M_g(x,\zeta_3,{x-1\over t})$ are of
order $(\log|g_1g_2|)x\log^{-2}x$, where the implied constant
is absolute. Furthermore, we obtain that
$M_g(x,\zeta_2,\zeta_3)=O({x\over \varphi(t)}{\log \log x\over \log^2 x})$, 
where again the implied constant is absolute. For the main term
we find that
$$M_g(x,\zeta_1)=A(g,t){x\over \log x}+
O\left({x\log|g_1g_2|\over \log^{2}x}\right),$$
where the implied constant is absolute.
On adding the various terms, the theorem follows. \qed\\

\indent The following result is a slight generalisation of
Lemma 2.4 of \cite{CM4}.
\begin{Lem} {\rm (GRH)}. 
\label{staartje}
Let $\psi(x)$ be a monotonous increasing positive
function which satisfies
$$\lim_{x\rightarrow \infty}\psi(x)=+\infty{\rm ~and~}\psi(x)\ll (\log x)^{1/2}.$$
Then we have 
$$\#\{p\le x: r_g(p)\ge \psi(x)\}\ll \log|g_1g_2|{\pi(x)\over \psi(x)},$$
where the constants implied by the $\ll$-symbol are absolute
and $g=g_1/g_2$ with $g_1$ and $g_2$ integers.
\end{Lem}
{\it Proof}. Let $y$ denote the largest integer not exceeding $\psi(x)$. We have
$$ \#\{p\le x: r_g(p)\ge y\}=\{p\le x:\nu_p(g)=0\}-\cup_{t=1}^{y-1}
R_g(0,1;t)(x),$$
where $\cup_{n=1}^{y-1}$ is a disjoint union. The latter
identity together with Theorem \ref{precies}, Lemma \ref{tail} 
and (\ref{landau}), yields
$$\#\{p\le x: r_g(p)\ge y\}=\pi(x)+O(\log|g_1g_2|)$$
$$-\left(1+O({h\over y})\right){x\over \log x}
+O\left({x\log y\log\log x\over \log^2 x}\right)+
O\left({xy\log|g_1g_2|\over\log^2 x}\right),$$
where the implied constants are all absolute. On using the prime number
theorem in the form $\pi(x)=x/\log x+O(x\log^{-3/2}x)$, the result then
follows. \qed
\subsection{On the convolution of the M\"obius function with Dirichlet
characters}
Let $\chi$ be a Dirichlet character of 
conductor $f_{\chi}$ and order
$o_{\chi}$ (for definitions and basic facts
on Dirichlet characters we refer the reader to Hasse \cite{Hasse0}).
Let $G_d$ denote the group of characters defined on $(\mathbb Z/d\mathbb Z)^*$.
We have that $G_d\cong (\mathbb Z/d\mathbb Z)^*$.
 An important auxiliary function in this paper is
the convolution, $h_{\chi}=\mu \star \chi$, of the M\"obius function  $\mu$ with
a Dirichlet character $\chi$, that is $h_{\chi}(n)=\sum_{d|n}\mu(d)\chi(n/d)$.
In this section we collect some auxiliary results involving $h_{\chi}$.
We note the following
trivial result.
\begin{Lem}
\label{vauxiliary}
The function $h_{\chi}$ is multiplicative. 
With the
convention that $0^{0}=1$, it satisfies $h_{\chi}(1)=1$
and $h_{\chi}(p^r)=\chi(p)^{r-1}[\chi(p)-1]$.
\end{Lem}
In particular if $\chi$ is the trivial character mod $d$, then
$$h_{\chi}(v)=\cases{\mu(v) &if $v|d$;\cr 0 &otherwise.}$$
By using one of the orthogonality relations for Dirichlet characters, the
following result is easily obtained.
\begin{Lem}
\label{nonmulttomult}
Let $a({\rm mod~}d)$ be a reduced residue class mod $d$. We have
$$\sum_{t\equiv a({\rm mod~}d)\atop t|v}\mu({v\over t})=
{1\over \varphi(d)}\sum_{\chi\in G_{d}}{\overline{\chi(a)}}h_{\chi}(v),$$
where $\chi_k$ runs over the Dirichlet characters modulo $d$.
\end{Lem}
Note that the lemma expresses a non-multiplicative function as a linear
combination of multiplicative functions. This will play an important
role later on.\\
\indent Let $r,s$ be non-negative integers. Put
$$C_{\chi}(h,r,s)=\sum_{(r,v)=1,~s|v}^{\infty}{h_{\chi}(v)(h,v)\over
v\varphi(v)}{\rm ~and~}A_{\chi}=\prod_{p\atop \chi(p)\ne
0}\left(1+{[\chi(p)-1]p\over [p^2-\chi(p)](p-1)}\right).$$ 
It is easy to see that the latter series is absolutely convergent.
Note that $h_{\chi}(v)\le 2^{\omega(v)}$, where $\omega(v)$ denotes the
number of distinct prime divisors of $v$. Note also that for every $\epsilon>0$
we have $|h_{\chi}(v)|\le 2^{\omega(v)}\le \sum_{d|v}1\ll v^{\epsilon}$ and
$\varphi(v)\gg v^{1-\epsilon}$. From this the absolute convergence
easily follows. Since $h_{\chi}$ is a multiplicative function, we can
invoke Euler's identity. After some tedious but easy calculations this then
yields the following result.
\begin{Lem}
\label{eulerproduct}
Let $h,r,s\ge 1$ be integers. Let $\chi$ be a Dirichlet character mod $d$.
Put $e_p=\nu_p(h)$.\\
i) If $(r,s)=1$, then
$$C_{\chi}(h,r,s)=\prod_{p\nmid
rs}\left(1+p^{1-e_p}{[\chi(p)-1]\over (p-1)}
\left({p^{e_p}-\chi(p^{e_p})\over p-\chi(p)}+
{\chi(p^{e_p})\over p^2-\chi(p)}\right)\right)$$
$$\prod_{\nu_p(s)=1}p^{1-e_p}{[\chi(p)-1]\over (p-1)}
\left({p^{e_p}-\chi(p^{e_p})\over p-\chi(p)}+
{\chi(p^{e_p})\over p^2-\chi(p)}\right)$$
$$\prod_{\nu_p(s)\ge e_p+1\atop \nu_p(s)\ge 2}\chi(p)^{\nu_p(s)-1}{[\chi(p)-1]\over (p-1)}{p^{e_p+3-2\nu_p(s)}\over 
(p^2-\chi(p))}$$
$$\prod_{2\le \nu_p(s)\le e_p}
{[\chi(p)-1]\over (p-1)}\left({(\chi(p)/p)^{\nu_2(s)-1}-(\chi(p)/p)^{e_p}\over 1-\chi(p)/p}
+{p^{1-e_p}\chi(p^{e_p})\over p^2-\chi(p)}\right).$$
In particular $C_{\chi}(h,r,s)=cA_{\chi}=cC_{\chi}(1,d,1)$, with 
$c\in \mathbb Q(\zeta_{o_{\chi}})$.\\
ii) If $(r,s)>1$, then $C_{\chi}(h,r,s)=0$.
\end{Lem}
{\tt Remark}. In this paper we only need to evaluate
$C_{\chi}(h,r,s)$ in the case where the largest odd divisor of $s$ is
squarefree and $\nu_2(s)\ge e_2+1$, in which case it assumes a bit
simpler form than the general one given in Lemma \ref{eulerproduct}.\\
 
\indent Only the primes $p$ with $\chi(p)\ne 1$ contribute to $A_{\chi}$.
Note that $A_{\chi}\ne 0$.
If $\chi$ is a principal character, then $A_{\chi}=1$.
If $\chi$ is real, then $A_{\chi}>0$.
If $\chi$ is a Dirichlet character and $\chi'$ is its associated primitive
Dirichlet character, then clearly $A_{\chi}$ equals $A_{\chi'}$ with some
local factors involving only $\zeta_{o_{\chi}}$ left out.
Hence
$A_{\chi}=c'A_{\chi'}$ with $c'\in \mathbb Q(\zeta_{o_{\chi}})$. Also note that $\overline A_{\chi}=A_{\overline 
\chi}$.\\
\indent The constants $A_{\chi}$ are the basic constants in this paper. They have
a product expansion in terms of special values of Dirichlet $L$-series \cite{Moreeaverage, PFibo}.
This result is related to the denominator identities arising in the theory of Lie algebras \cite{Moreewitt}.
These expansions can be used to evaluate $A_{\chi}$ with high numerical accuracy \cite{Moreeaverage}. The
values of $A_{\psi_1}$ and $A_{\xi_1}$ given in the introduction are taken from Table 3 of \cite{Moreeaverage}.\\
\indent Another result involving $h_{\chi}$ needed is the following (where the sums are over the integers $v\ge 1$).
\begin{Lem}
\label{laatsteh}
Let $r,s$ be integers with $s|r$ and $g>0$. Let $\chi$ be a 
Dirichlet character. Then, if $g>0$ or $g<0$ and $s$ is even,
$$\sum_{(r,v)=1}{h_{\chi}(v)\over [K_{sv,v}:\mathbb Q]}=
{1\over 
\varphi(s)}\left(C_{\chi}(h,r,1)+C_{\chi}(h,r,{n_s\over (n_s,s)})\right).$$
When $g<0$ and $s$ is odd, the latter sum equals
$$
{1\over 
\varphi(s)}\left(C_{\chi}(h,r,1)-{1\over 2}C_{\chi}(h,r,2)
+{1\over 2}C_{\chi}(h,r,2^{\nu_2(h)+1})+
C_{\chi}(h,r,{n_s\over (n_s,s)})\right).$$
\end{Lem}
{\it Proof}. For $g>0$ or
$g<0$ and $2|s$, the
 proof easily follows from the identity
$$\sum_{(r,v)=1}{h_{\chi}(v)\over [K_{sv,v}:\mathbb Q]}=
\sum_{(r,v)=1}{h_{\chi}(v)(h,v)\over v\varphi(sv)}+
\sum_{(r,v)=1\atop n_s|sv}{h_{\chi}(v)(h,v)\over v\varphi(sv)},$$
which on its turn is an easy consequence of Lemma \ref{degree}. 
The proof of the remaining case is similar. \qed

\subsection{Preliminaries specific to the case $d=4$}
\subsubsection{Even order}
\label{deordeiseven}
Let $s\ge 2$.
It is not difficult to estimate $N_g(1,2^s;0,4)(x)$. To this end we
consider the set of primes $p\equiv 1({\rm mod~}2^s)$ such that
ord$_g(p)\not\equiv 0({\rm mod~}4)$. Let $r=\nu_2(p-1)$. Note that 
$r\ge s$. Now $p$
satisfies ord$_g(p)\not\equiv 0({\rm mod~}4)$ if and
only if
$g^{(p-1)/2^{r-1}}\equiv 1({\rm mod~}p)$, that is if and only if
$p$ splits completely in $K_{2^r,2^{r-1}}$, but not completely in
$K_{2^{r+1},2^{r-1}}$. By Chebotarev's density theorem we expect that
$N_g(1,2^s;0,4)$ has a density as given in (\ref{dens0mod4}).
Indeed, this can be shown unconditionally. Wiertelak's work
\cite[Theorem 2]{Wiertelak1} goes beyond this and shows that
we even have an estimate for $N_g(1,2^s;0,4)(x)$ as given by
(\ref{densprecise}).\\
\indent Similarly we can easily estimate
$N_g(1,2^s;1,4)(x)+N_g(1,2^s;3,4)(x)$ which equals $N_g(1,2^s;1,2)(x)$ and
$N_g(1,2^s;2,4)(x)$. See Theorem \ref{main1} for the outcome.

\subsubsection{Odd order}
We let $\psi_0,\psi_1$ denote the trivial
respectively non-trivial character mod $4$.
The starting point of our analysis is the following easy
result.
\begin{Lem}
\label{startingpoint}
Let $s\ge 1$. For $j$ odd we have
$$N_g(1,2^s;j,4)(x)=\sum_{r\ge s}\sum_{t\equiv j({\rm mod~}4)}
\#\{p\le x:p\equiv 1+2^r({\rm mod~}2^{r+2}),~r_g(p)=2^rt\}$$
$$+\sum_{r\ge s}\sum_{t\equiv -j({\rm mod~}4)}
\#\{p\le x:p\equiv 1+3\cdot 2^r({\rm mod~}2^{r+2}),~r_g(p)=2^rt\},$$
$$N_g(3,8;j,4)(x)=\sum_{t\equiv j({\rm mod~}4)}\#\{p\le x:p\equiv 3({\rm
mod~}8),~r_g(p)=2t\}{\rm ~and~}$$$$N_g(7,8;j,4)(x)=\sum_{t\equiv -j({\rm
mod~}4)}\#\{p\le x:p\equiv 7({\rm
mod~}8),~r_g(p)=2t\}.$$\end{Lem}
{\it Proof}. We only prove the assertion regarding $N_g(1,2^s;j,4)(x)$,
the other assertions being easier to prove. For every prime
$p\equiv 1({\rm mod~}2^s)$ there exists an unique $r\ge s$ such that
either $p\equiv 1+2^r({\rm mod~}2^{r+2})$ or
$p\equiv 1+3\cdot 2^r({\rm mod~}2^{r+2})$. We assume that we are
in the first case, the other case being dealt with similarly.
Using (\ref{indord}) we note that ord$_g(p)\equiv j({\rm mod~}4)$
if and only if $r_g(p)=2^rt$ for some $t\ge 1$ with $t\equiv j({\rm mod~}4)$.
\qed\\

\noindent The densities of the
sets appearing in Lemma \ref{startingpoint} can be
determined, under GRH, on invoking (\ref{lenstraoud}). Assuming the
densities add up, we then arrive at the conjecture that the densities are as stated
in Theorem \ref{triplesum} (cf. Section \ref{proofofthemainresults}).\\
\indent An alternative approach is obtained on first resumming the
expressions in Lemma \ref{startingpoint} and only then applying
Lenstra's machinery \cite{Lenstra}, which is what is explained next.\\
\indent Let $s\ge 1$. From Lemma \ref{startingpoint} and (\ref{inclusionex}) we
deduce that, when $j$ is odd, $$N_g(1,2^s;j,4)(x)=
\sum_{r\ge s}\sum_{t\equiv j({\rm mod~}4)}\sum_{n=1}^{\infty}\mu(n)
S_g(1+2^r,2^{r+2};nt2^r,nt2^r)(x)$$
\begin{equation}
\label{lang}
+\sum_{r\ge s}\sum_{t\equiv -j({\rm mod~}4)}\sum_{n=1}^{\infty}\mu(n)
S_g(1+3\cdot 2^r,2^{r+2};nt2^r,nt2^r)(x).
\end{equation}
Note that the fields that arise in the $S_g$
occurring in (\ref{lang}) are of
the form $K_{v2^r,v2^r}$. On grouping together the contributions involving
the various $K_{v2^r,v2^r}$ the triple sums can be reduced
to double sums. To this end
we first note that we can restrict to the case where $n$ is odd, since
if $n$ is even and $m$ is
odd, $S_g(1+m\cdot 2^r,2^{r+2};nt2^r,nt2^r)$ is empty (then the
condition on the Frobenius symbol implies that $p\equiv 1({\rm
mod~}2^{r+1})$). On putting $v=nt$ the summation is then over 
all odd $v\ge 1$.
As weighing factors we then get sums as in Lemma \ref{nonmulttomult} with
$(a,f)=(1,4)$ and $(a,f)=(3,4)$. On applying Lemma \ref{nonmulttomult} and
noting that for odd $v$, $h_{\psi_0}(v)=1$ if $v=1$ and
$h_{\psi_0}(v)=0$ otherwise, we then obtain that
$N_g(1,2^s;j,4)(x)={1\over 2}I_1+{(-1)^{j-1\over 2}\over 2}I_2$, where
$$I_1=\sum_{r\ge s}\Big[\#\{p\le x: (p,K_{2^r,2^r}/\mathbb Q)=id\}
- \#\{p\le x: (p,K_{2^{r+1},2^r}/\mathbb Q)=id\}\Big]$$
and
$$I_2=\sum_{r\ge s}\sum_{2\nmid v}^{\infty}h_{\psi_1}(v)
\Big[S_g(1+2^r,2^{r+2};v2^r,v2^r)(x)-S_g(1+3\cdot
2^r,2^{r+2};v2^r,v2^r)(x)\Big].$$ 
(Note that the latter double sum can be simplified
to a single sum. On doing so we find that $I_2$ equals
$$\sum_{2^s|w}h_{\psi_1}(w_{odd})
\Big[S_g(1+2^{\nu_2(w)},2^{\nu_2(w)+2};w,w)(x)-
S_g(1+3\cdot 2^{\nu_2(w)},2^{\nu_2(w)+2};w,w)(x)\Big],$$
where $w_{odd}$ is the largest odd divisor of $w$.)
Note that if we add
$N_g(1,2^s;1,4)(x)$ and $N_g(1,2^s;3,4)(x)$ we obtain $I_1$,
which is a well-known result.\\
\indent If $s\ge 2$, then the Chebotarev density theorem implies,
unconditionally, that for $v$ odd and $r\ge s$,
\begin{equation}
\label{hetzelfde}
S_g(1+2^r,2^{r+2};v2^r,v2^r)(x)\sim S_g(1+3\cdot
2^r,2^{r+2};v2^r,v2^r)(x),{\rm ~as~}x\rightarrow \infty.
\end{equation}
Thus we might expect that $I_2$ behaves
like an error term and that, consequently,
$N_g(1,2^s;1,4)(x)\sim I_1/2$ as $x$ tends to infinity. Theorem
\ref{main1} shows that this is indeed true, under GRH.\\
\indent If $s=1$, however, then (\ref{hetzelfde}) does not necessarily
 hold true for every $r\ge s$.
It thus makes sense to consider $N_g(3,4;j,4)(x)$ for $j$ odd
separately. We then
obtain analogous expressions to those for $I_1$ and
$I_2$, but instead of summing over $r\ge s$ we take $r$ to equal one:
$$N_g(3,4;j,4)(x)=
{1\over 2}\#\{p\le x:p\equiv 3({\rm mod~}4),~({g\over p})=1\}+$$
$${(-1)^{j-1\over 2}\over 2}\sum_{2\nmid v}^{\infty}h_{\psi_1}(v)
\Big[S_g(3,8;2v,2v)(x)-S_g(7,8;2v,2v)(x)\Big].$$
Using that $p$ splits completely in $\mathbb Q(\sqrt{-2})$ iff
$p\equiv 1({\rm mod~}8)$ or $p\equiv 3({\rm mod~}8)$ and that
$p$ splits completely in $\mathbb Q(\sqrt{2})$ iff
$p\equiv \pm 1({\rm mod~}8)$, we obtain the following result.
We present it with a more succinct proof. For a number field $L$
we let $\pi_L(x)$ denote the number of rational primes $p\le x$
that split completely in $L$. 
\begin{Lem}
\label{simpeltoch}
Let $j$ be odd. For every $x$ we have
$$
N_g(3,4;j,4)(x)={1\over 2}\#\{p\le x:p\equiv 3({\rm mod~}4),({g\over p})=1\}+
{(-1)^{j-1\over 2}\over 2}\sum_{2\nmid v}^{\infty}h_{\psi_1}(v)A_v(x),
$$
where
$
A_v(x)=\pi_{K_{2v,2v}(\sqrt{-2})}(x)-\pi_{K_{2v,2v}(\sqrt{2})}(x)$.
\end{Lem}
{\it Proof}. 
Let us consider only the case where $j\equiv 1({\rm mod~}4)$, the remaining
case being dealt with similarly.
Note that in both the left hand side and the right
hand side of the identity that is to be established only primes $p$ satisfying
$p\equiv 3({\rm mod~}4)$ and $\nu_p(g)=0$ are counted. Now let $p$
be a prime such that $p\equiv 3({\rm mod~}4)$, $\nu_p(g)=0$ and $p\le x$.
We shall show that it is counted with the same
multiplicity in both the left and the right hand side of the identity that
is to be established,
thus
finishing the proof.\\
\indent Since by assumption $\nu_p(g)=0$, there exists a largest integer
$k$ such that $g^{p-1\over k}\equiv 1({\rm mod~}p)$. Note that in both
the left and the right hand side only primes $p$ with $k$ even are
counted. Thus we may write $k=2k_1$. Note that $k_1$ must be odd.
Let us assume that $p\equiv 3({\rm mod~}8)$. Then $p$ is counted on the
right hand side with weight
$${1\over 2}+{1\over 2}\sum_{m|k_1}h_{\psi_1}(m)
={1\over 2}+{1\over 2}(\psi_1\star \mu \star {\bf
1})(k_1)={1+\psi_1(k_1)\over 2}.$$
Thus the weight is $1$ if $k_1\equiv 1({\rm
mod~}4)$ and $0$ otherwise. In
other words
the weight is $1$ iff ord$_g(p)\equiv 1({\rm mod~}4)$.\\
\indent The case
where  $p\equiv 7({\rm mod~}8)$ is dealt with similarly. \qed\\

Using Chebotarev's density theorem we now expect that
$$\lim_{x\rightarrow \infty}{\sum_{2\nmid
v}^{\infty}h_{\psi_1}(v)A_v(x)\over \pi(x)}=\sum_{2\nmid
v}\Big({h_{\psi_1}(v)\over [K_{2v,2v}(\sqrt{-2}):\mathbb
Q]}-{h_{\psi_1}(v)\over [K_{2v,2v}(\sqrt{2}):\mathbb Q]}\Big).$$ Let us denote the
quantity on the right hand side by
${\tilde \Delta_g}$. Note that
$${\tilde \Delta_g}={1\over 2}\sum_{\sqrt{-2}\in K_{2v,2v}\atop 2\nmid
v}{h_{\psi_1}(v)\over [K_{2v,2v}:\mathbb Q]}-{1\over 2}\sum_{\sqrt{2}\in
K_{2v,2v}\atop 2\nmid v}{h_{\psi_1}(v)\over
[K_{2v,2v}:\mathbb Q]}={\Delta_g\over 2}.$$
Theorem \ref{main3} shows that this heuristic holds true, under GRH.\\

\section{The distribution of the index over residue classes}
\label{indexmodd}
The problem of the distribution of the index over residue classes is
far easier than that of the distribution of the order. However, the
answers to both problems turn out to have some features in common.\\
\indent Let $a$ and $d$ be integers. Under GRH it follows from
Pappalardi's work \cite{Pappalardi} that
the density, $\rho_g(a,d)$, of the set of primes $p$ such that
$r_g(p)\equiv a({\rm mod~}d)$ exists and equals
$$\rho_g(a,d)=\sum_{t\equiv a({\rm mod~}d)}\sum_{v=1}^{\infty}
{\mu(v)\over [K_{vt,vt}:\mathbb Q]}.$$
Using this the following result is then easily deduced.
\begin{Thm} {\rm (GRH)}.
\label{indeks}
Let $a$ and $d$ be arbitrary natural numbers. Put $\delta=d/(a,d)$. Then
the density of the primes $p$ 
with $r_g(p)\equiv a({\rm mod~}d)$, $\rho_g(a,d)$,  exists and satisfies
$$\rho_g(a,d)=\sum_{\chi\in G_{\delta}}c_{\chi}A_{\chi},~{\rm ~with~}c_{\chi}
\in \mathbb Q(\zeta_{o_{\chi}}).$$
Furthermore, $c_{\overline \chi}={\overline c_{\chi}}$. The number
$c_{\chi}$ can be explicitly computed.
\end{Thm}
{\it Proof}. On putting $w=(a,d)$ and $\alpha=a/w$ we obtain
$$\rho_g(a,d)=\sum_{t\equiv \alpha({\rm mod~}\delta)}
\sum_{v=1}^{\infty}{\mu(v)\over [K_{vwt,vwt}:\mathbb Q]}.$$
Writing $vt=v_1$ and invoking Lemma \ref{nonmulttomult}, we then obtain
$$\rho_g(a,d)={1\over \varphi(\delta)}\sum_{\chi \in G_{\delta}}
{\overline{\chi(\alpha)}}\sum_{v_1=1}^{\infty}{h_{\chi}(v_1)\over
[K_{v_1w,v_1w}:\mathbb Q]}.$$
Let $g>0$. By Lemma \ref{degree} we have
$$\sum_{v=1}^{\infty}{h_{\chi}(v)\over [K_{vw,vw}:\mathbb Q]}
=\sum_{v=1}^{\infty}{h_{\chi}(v)(h,vw)\over vw\varphi(vw)}+
\sum_{v=1\atop n_1|vw}^{\infty}{h_{\chi}(v)(h,vw)\over vw\varphi(vw)}
=J_1+J_2.$$
We can rewrite $J_1$ as
$$J_1={(h,w)\over w\varphi(w)}\sum_{v=1}^{\infty}{h_{\chi}(v)
(h,vw)\varphi(w)\over (h,w)v\varphi(vw)},$$
where the argument of the sum is easily seen to be multiplicative in $v$.
If $p\nmid hwf_{\chi}$, then the local factor at $p$ in the Euler
product for $J_1$ equals that of $A_{\chi}$ and so $J_1=c'_{\chi}A_{\chi}$, 
where $c'_{\chi}\in \mathbb Q(\zeta_{o_{\chi}})$. Rewriting the condition
$n_1|vw$ as ${n_1\over (n_1,w)}|v$, we see that also $J_2=c''_{\chi}A_{\chi}$, 
where
$c''_{\chi}\in \mathbb Q(\zeta_{o_{\chi}})$. 
A similar argument can be used in the case where $g<0$, cf. the
proof of Lemma \ref{laatsteh}.\\
\indent On using that ${\overline h_{\chi}}=h_{\overline \chi}$ and
${\overline A_{\chi}}=A_{\overline \chi}$ it follows that $c_{\overline \chi}={\overline c_{\chi}}$. \qed\\

\noindent A special case occurs when $d|a$. Then $\rho_g(a,d)$ is the
density of primes $p\le x$ such that the index of $g({\rm mod~}p)$ is
divisible by $d$, that is $\rho_g(a,d)$ is the density of primes $p\le x$
that split completely in $K_{d,d}$. By an unconditional version of the
Chebotarev density theorem we then infer
\begin{Prop}
\label{eenvoudigzeg}
The density $\rho_g(0,d)$ exists and satisfies
$$\rho_g(0,d)={1\over [K_{d,d}:\mathbb Q]}.$$
\end{Prop}
\indent In the three examples some special cases of Theorem \ref{indeks} are discussed.\\ 
\noindent {\bf Example 1}.  (GRH). We consider the case where $(a,d)=1$ and $g>0$.
Then, using Lemma \ref{laatsteh}, we obtain that
$$\rho_g(a,d)={1\over \varphi(d)}\sum_{\chi\in G_d}{\overline{\chi(a)}}
\Big(C_{\chi}(h,1,1)+C_{\chi}(h,1,n_1)\Big).$$

\noindent {\bf Example 2}. (GRH). We assume that $(a,d)=1$ and $g=2$. Note 
that $D=8$ and
hence $n_1=8$. On invoking the formula of Example 1 with $h=1$ and $n_1=8$, we
obtain, using Lemma \ref{eulerproduct},
$$\rho_2(a,d)={1\over \varphi(d)}\sum_{\chi \in G_d}{\overline 
{\chi(a)}}\cdot \delta_{\chi},$$
with
$$\delta_{\chi}=\left({1\over 2}+
{\chi_(2)((\chi(2))^2-\chi(2)+12)\over 8(4-\chi(2))}\right)
\prod_{p>2}\left(1+{p(\chi(p)-1)\over (p-1)(p^2-\chi(p))}\right).$$
(This corrects a typo in Corollary 8 of \cite{Pappalardi}). Alternatively we can write
$$\rho_2(a,d)={1\over \varphi(d)}\prod_{p|d}\left(1-{1\over p(p-1)}\right)\sum_{\chi \in G_d}{\overline 
{\chi(a)}}\left(1+{\chi(4)(\chi(2)-1)\over 8(\chi(2)+2)}\right)A_{\chi}.$$

\noindent {\bf Example 3} (GRH). Let $g=2$ and $d=3$. Using Example 2 we 
compute $\rho_2(\pm 1,3)={5\over 12}\pm {5\over 16}A_{\xi_1}$,
that is $\rho_2(1,3)=0.471034\cdots$ and $\rho_2(2,3)=0.362298\cdots$.
(This corrects the values claimed in \cite[p. 386]{Pappalardi}.) 
From the latter two values we infer that $\delta_2(0,3)=1/6$.
By Proposition \ref{eenvoudigzeg} it follows that we even have
unconditionally that $\delta_2(0,3)=1/6$.
Up to
$p_{10^6}=1299709$ the approximations $0.16589$, $0.47127$ and
$0.36283$ for
$\rho_2(0,3)$, $\rho_2(1,3)$ respectively $\rho_2(2,3)$ are found.

\section{Proof of the main results}
\label{proofofthemainresults}
In our formulation of the main results stated in the introduction, we
have added the (known) case $d=2$ for completeness. As 
already pointed
out in Section \ref{deordeiseven} these cases have been well studied.
The best known error terms are due to Wiertelak, cf.
\cite{Wiertelak2}. The densities for these cases can be explicitly 
evaluated using Lemma \ref{degree}. Since the interested reader can
easily carry this out herself, we abstained from writing down the rather
lengthy (because of case distinctions) outcome.
For some further elaboration on these cases see
Section 2.5.1.\\
\indent Our proof for the remaining cases has an 
analytic and algebraic component, with the analytic component 
being captured by the
following result. 

\begin{Thm}
\label{triplesum}
{\rm (GRH)}. Let $s\ge 1$. 
Let $\psi_1$ be the non-principal character modulo $4$. For $r\ge 1$
let $\sigma_{1,r},\sigma_{-1,r}
\in {\rm Gal}(\mathbb Q(\zeta_{2^{r+2}})/\mathbb Q)$ be the automorphisms
that are uniquely determined by
$\sigma_{1,r}(\zeta_{2^{r+2}})=\zeta_{2^{r+2}}^{1+2^r}$, respectively
$\sigma_{-1,r}(\zeta_{2^{r+2}})=\zeta_{2^{r+2}}^{1+3\cdot 2^r}$.
For $j=-1$ and $j=1$ let 
$$
c_j(r,tn)=
\cases{1 &if $\sigma_{j,r}|_{\mathbb Q( \zeta_{2^{r+2}})\cap \mathbb
Q(\zeta_{2^rtn},g^{1/2^rtn})}=$id;\cr
0 &otherwise.}
$$
For $j=1$ and $j=3$ we have
$$N_g(1,2^s;j,4)(x)=\delta_g(1,2^s;j,4){x\over \log x}+
O\left(\log|g_1g_2|{x\over \log^{3/2}x}\right),$$
where
$$\delta_g(1,2^s;j,4)=\sum_{r\ge s}^{\infty}\sum_{t=1\atop 2\nmid t}^{\infty}
\sum_{n=1}^{\infty}{\mu(n)c_{\psi_1(jt)}(r,tn)
\over [\mathbb Q(\zeta_{2^{r+2}},\zeta_{2^rtn},g^{1/2^rtn}):\mathbb Q]},$$
and the implied constant is absolute.
\end{Thm}
An heuristic argument in favour of the truth of the latter theorem is
easily given. 
From Lemma \ref{startingpoint} and (\ref{inclusionex}) we
deduce that $N_g(1,2^s;j,4)(x)$ equals
$$\sum_{r\ge s}\sum_{t\equiv j({\rm mod~}4)}\sum_{n=1}^{\infty}\mu(n)
\#\{p\le x:p\equiv 1+2^r({\rm mod~}2^{r+2}),~(p,K_{2^rtn,2^rtn}/\mathbb
Q)=id\}+$$ 
$$\sum_{r\ge s}\sum_{t\equiv -j({\rm mod~}4)}\sum_{n=1}^{\infty}\mu(n)
\#\{p\le x:p\equiv 1+3\cdot 2^r({\rm
mod~}2^{r+2}),~(p,K_{2^rtn,2^rtn}/\mathbb Q)=id\}.$$
The density of the inner
sums is given, under GRH, by (\ref{lenstraoud}).
Assuming the densities add up and there
is sufficiently cancellation in the error terms, we then arrive at 
the heuristic that the densities should be as claimed. Our proof
of Theorem \ref{triplesum} is in the same spirit:\\

\noindent {\it Proof of Theorem} \ref{triplesum}. Let us denote the triple sum in the formulation of the result by
$\sum \sum \sum$. All constants implied by the O-symbols in this proof
will be absolute.
The first formula of Lemma \ref{startingpoint} can be more compactly
written as
$$N_g(1,2^s;j,4)(x)=\sum_{r\ge s}
\sum_{2\nmid t}R_g(1+(2-\psi_1(jt))2^r,2^{r+2},2^rt)(x).$$
On retaining only the primes with $r_g(p)\le y$, we obtain
\begin{equation}
\label{ophakken}
N_g(1,2^s;j,4)(x)=T_1(y)+O(\#\{p\le x:r_g(p)\ge y\}),
\end{equation}
where 
$$T_1(y)=\sum_{r\ge s}\sum_{2\nmid t\atop 2^rt\le y}R_g(1+(2-\psi_1(jt))2^r,2^{r+2},2^rt)(x).
$$
The function $R_g(1+(2-\psi_1(jt))2^r,2^{r+2},2^rt)(x)$ can be estimated
in the same way
as $R_g(0,1;t)(x)$ (see the proof of Theorem \ref{precies}). We find
that, under GRH,
$$R_g(1+(2-\psi_1(jt))2^r,2^{r+2},2^rt)(x)={\rm Li}(x)
\sum_{n=1}^{\infty}{\mu(n)c_{\psi_1(jt)}(r,tn)\over 
[\mathbb Q(\zeta_{2^{r+2}},\zeta_{2^rtn},g^{1/2^rtn}):\mathbb Q]}$$
$$+O\left(\log(|g_1g_2|){x\over \log^2 x}\right)
+O\left({x\log\log x\over \varphi(2^rt)\log^{2}x}\right).
$$
This, when substituted in (\ref{ophakken}), yields 
$$T_1(y)={\rm Li}(x)\sum\sum\sum + O\left({h{\rm Li}(x)\over y}\right)
+O\left({x\log y\log\log x\over \log^2 x}\right)+
O\left({xy\log|g_1g_2|\over\log^2 x}\right),$$
where we used (\ref{landau}) and
$$
\sum_{r\ge s}\sum_{2\nmid t\atop 2^rt>y}\sum_{n=1}^{\infty}
{h\over 2^rtn\varphi(2^rtn)}=
O\left(\sum_{m\ge y}\sum_{n=1}^{\infty}{h\over mn\varphi(mn)}\right)
=O({h\over y}),$$
cf. the proof of Lemma \ref{tail}. On taking $y=\sqrt{\log x}$
in (\ref{ophakken}), the result then follows on invoking
Lemma \ref{staartje} with $\psi(x)=\sqrt{\log x}$.\qed\\ 

The algebraic part is a consequence of the following result.
\begin{Lem}
\label{bijnagelijk}
Let $n$ be squarefree and $t$ be odd. Denote the intersection of
$\mathbb Q(\zeta_8)$ and $K_{2nt,2nt}$ by $L_{nt}$.\\
{\rm i)} If $r=1$, $2\nmid hn$, $8|D$ and $D|8nt$, then
$$L_{nt}=\mathbb Q(\sqrt{2\cdot {\rm sgn}(g)}~),~c_1(r,tn)={1-{\rm sgn}(g)
\over 2},~c_{-1}(r,tn)={1+{\rm sgn}(g)\over 2},$$ 
if $D\equiv 8({\rm mod~}32)$ and
$$L_{nt}=\mathbb Q(\sqrt{2\cdot {\rm sgn}(-g)}~),~c_1(r,tn)=
{1+{\rm sgn}(g)\over 2},~c_{-1}(r,tn)={1-{\rm sgn}(g)\over 2} ,$$ 
if $D\equiv 24({\rm mod~}32)$.\\
{\rm ii)} We have $c_1(r,tn)\ne c_{-1}(r,tn)$ if and 
only if $r=1$, $2\nmid hn$, $8|D$ and $D|8nt$. 
\end{Lem}
{\it Proof}. i). Under the hypothesis of part i) we can apply Lemma
\ref{intersection} to infer that $L_{nt}$, the field intersection in the
definition of $c_{\pm 1}(1,tn)$, equals $\mathbb Q(\sqrt{2\cdot {\rm sgn}(g)}~)$ if
$D\equiv 8({\rm mod~}32)$ and  $\mathbb Q(\sqrt{2\cdot {\rm sgn}(-g)}~)$ in the
remaining case $D\equiv 24({\rm mod~}32)$. Writing $\sqrt{2}=\zeta_8+\zeta_8^{-1}$ one
calculates that $\sigma_{\pm 1,1}(\sqrt{2})=\mp \sqrt{2}$,
$\sigma_{\pm 1,1}(\sqrt{-2})=\pm \sqrt{-2}$. From this observation
and the definition of $c_{\pm 1}(r,tn)$ the result follows at once.\\
ii). `$\Leftarrow$'. Follows by
part i). `$\Rightarrow$'. The intersection of the fields
$\mathbb Q( \zeta_{2^{r+2}})$ and
$\mathbb Q(\zeta_{2^rtn},g^{1/2^rtn})$ is abelian and, since
$4\nmid nt$, is contained in
$\mathbb Q(\zeta_{2^{r+1}},\sqrt{2})$
or $\mathbb Q(\zeta_{2^{r+1}},\sqrt{-2})$.
As $\sigma_{1,r}(\zeta_{2^{r+1}})=\sigma_{-1,r}(\zeta_{2^{r+1}})$ for every $r\ge 1$, we deduce that
$\mathbb Q( \zeta_{2^{r+2}})\cap \mathbb Q(\zeta_{2^rtn},g^{1/2^rtn})$ must
contain
at least one element from $\{\sqrt{-2},\sqrt{2}\}$.\\
\indent Now let us consider how $\sigma_{1,r}$
and $\sigma_{-1,r}$ act on $\sqrt{2}~(=\zeta_8+\zeta_8^{-1})$. For $r\ge 2$
we have
$\sigma_{1,r}(\sqrt{\pm 2})=\sigma_{-1,r}(\sqrt{\pm 2})$, since then
$\sigma_{1,r}(\zeta_8)=
\sigma_{-1,r}(\zeta_8)$. 
Thus we must have $r=1$. If $n$ is even, then
$i\in \mathbb Q( \zeta_{2^{r+2}})\cap \mathbb Q(\zeta_{2^rtn},g^{1/2^rtn})$
and
since $\sigma_{\pm 1,1}(i)=-i$, we infer that
$c_{\pm 1}(r,tn)=0$. Thus $n$ is odd.\\
\indent From the above discussion it follows that 
$\mathbb Q(\zeta_8)\cap \mathbb Q(\zeta_{2nt},g^{1/2nt})$ 
with $nt$ odd must contain $\sqrt{2}$ or $\sqrt{-2}$.
By Lemma \ref{intersection} this leads then to the further restrictions
(apart from $r=1$ and $2\nmid n$): $2\nmid h$, $8|D$ and $D|8nt$. \qed\\

\noindent {\it Proof of Theorems {\rm \ref{main1}} and
{\rm \ref{main3}}}. 
The claims in the cases where the order is even or
$1({\rm mod~}2)$ are already known. 
The results for the remaining cases are a straigthforward consequence of
Theorem \ref{triplesum} and Lemma \ref{bijnagelijk}. If $s\ge 2$, then
$r\ge 2$ in the triple sum for the densities. By Lemma \ref{bijnagelijk}
we then have that $c_{-1}(t,rn)=c_{1}(t,rn)$ and hence
$\delta_g(1,2^s;1,4)=\delta_g(1,2^s;3,4)$. On noting that
$N_g(1,2^s;1,4)(x)+N_g(1,2^s;3,4)=N_g(1,2^s;1,2)(x)$ the proof of
Theorem \ref{main1} is then completed.\\
\indent Since for `most' $t$ and $n$ we have $c_1(1,tn)=c_{-1}(1,tn)$
(Lemma \ref{bijnagelijk}) it is natural to
compute the difference $\delta_g(3,4;1,4)-\delta_g(3,4;3,4)$. 
Since
$\delta_g(3,4;1,4)+\delta_g(3,4;3,4)$ is easily evaluated, we are then done. 
We proceed by filling in the details.\\
\indent Proceeding as in the proof of Theorem \ref{triplesum} we infer that
$$\delta_g(3,4;j,4)=\delta_g(1,2;j,4)-\delta_g(1,4;j,4)=
\sum_{t=1\atop 2\nmid t}^{\infty}\sum_{n=1}^{\infty}
{\mu(n)c_{\psi_1(jt)}(1,tn)\over
[\mathbb Q(\zeta_8,\zeta_{2tn},g^{1/2tn}):\mathbb Q]}$$
(now only the terms with $r=1$ contribute and hence the triple sum is reduced
to a double sum).
From the latter formula and the fact that $c_1(1,tn)=c_{-1}(1,tn)$
for $t$ is odd and $n$ is even (see Lemma \ref{bijnagelijk}), we 
infer that
$$
\delta_g(3,4;1,4)-\delta_g(3,4;3,4)=
\sum_{2\nmid t}\sum_{2\nmid n}{\mu(n)[c_{\psi_1(t)}(1,tn)-c_{-\psi_1(t)}(1,tn)]
\over
[\mathbb Q(\zeta_8,\zeta_{2tn},g^{1/2tn}):\mathbb Q]}.$$
On writing $nt=v$ we obtain
$$\delta_g(3,4;1,4)-\delta_g(3,4;3,4)=
\sum_{2\nmid v}{h_{\psi_1}(v)[c_1(1,v)-c_{-1}(1,v)]\over 
[\mathbb Q(\zeta_8,\zeta_{2v},g^{1/2v}):\mathbb Q]}.$$
On invoking Lemma \ref{bijnagelijk} the latter sum is seen
to equal $\Delta_g/2$. We 
now infer from Theorem \ref{triplesum} that
$$N_g(3,4;1,4)(x)-N_g(3,4;3,4)(x)={\Delta_g\over 2}{\rm Li}(x)
+O\left(\log|g_1g_2|{x\over \log^{3/2} x}\right),$$
where the implied constant is absolute.
On noting that
$$N_g(3,4;1,4)(x)+N_g(3,4;3,4)(x)=\#\{p\le x:p\equiv 3({\rm mod~}4),~
({g\over p})=1\},$$
the result easily follows. \qed

\section{Explicit evaluation of the densities}
\label{Explicit}
In this section we explicitly evaluate the
densities (under GRH), computed in Theorem \ref{main1} and
Theorem \ref{main3}.\\
\indent Using Lemma \ref{degree}, Theorem \ref{main1} can be made
completely explicit. For reasons of space we restrict ourselves to
describing the situation for a `generic' $g$.
\begin{Thm} {\rm (GRH)}. Let $s\ge 2$. If $h$ is odd and $D$ contains an
odd prime factor, then
$$\delta_g(1,2^s;j,4)=\cases{2^{1-s}-{2\over 3}\cdot 4^{1-s} &if $j=0$;\cr
4^{1-s}/6 &if $j=1$;\cr 4^{1-s}/3 &if $j=2$;\cr
4^{1-s}/6 &if $j=3$.}$$
\end{Thm}
{\it Proof}. The conditions on $h$ and $D$ ensure that the degrees
$[K_{m,n}:\mathbb Q]$ occurring in the sums in Theorem \ref{main1}
equal $\varphi(m)n$. It then remains to sum some geometric series. \qed\\

\indent Theorem \ref{main3} shows that the sets $N_g(3,4;j,4)$ considered
there have a density, $\delta_g(3,4;j,4)$, under GRH. The case
where $j$ is even is trivial and left to the reader.
\begin{Thm} {\rm (GRH)}.
\label{main3explicit}
Let $g\in \mathbb Q\backslash\{-1,0,1\}$, Write $g=\pm g_0^h$, where $g_0>0$ is not a 
power of
a rational number.
For any prime $p$ define $e_p$ by $p^{e_p}||h$.\\ 
\indent If $h$ is even, then $\delta_g(3,4;j,4)=(1+{\rm sgn}(g))/8$.\\
\noindent Next, let $h$ be odd.
Then
$\delta_g(3,4;1,4)=\delta_g(3,4;3,4)=1/8$, unless $D$, the discriminant of
the quadratic field $\mathbb Q(\sqrt{g_0})$, is divisible by $8$ and has no
prime divisor congruent to
$1({\rm mod~}4)$, in which case we have
$$\delta_g(3,4;j,4)={1\over 8}+{\rm sgn}(g){(-1)^{j+1\over 2}\over 8}P_1P_2P_3,$$
where
$$P_1=\prod_{p|D\atop p\equiv 3({\rm
mod~}4)}\left({2[p^{e_p}-(-1)^{e_p}]\over p^{e_p-1}(p^2-1)}
+{2p(-1)^{e_p}\over p^{e_p}(p^2+1)(p-1)}\right),
$$
$$P_2=\prod_{p\nmid D,~p|h\atop p\equiv 3({\rm mod~}4)}
\left(1-{2[p^{e_p}-(-1)^{e_p}]\over p^{e_p-1}(p^2-1)}
+{2p(-1)^{e_p+1}\over p^{e_p}(p^2+1)(p-1)}\right),~{and~}$$
$$P_3=\prod_{p\nmid hD\atop p\equiv 3({\rm mod~}4)}
\left(1-{2p\over (p^2+1)(p-1)}\right),$$
and thus in particular {\rm sgn}$(\delta_g(3,4;3,4)-\delta_g(3,4;1,4))=
{\rm sgn}(g)$, since the local factors of $P_1P_2P_3$ are all positive.
\end{Thm}
\begin{Cor} {\rm (GRH)}.
For $j=0,1,2,3$ we have
$\delta_g(3,4;j,4)=c_1(j)+c_2(j)A_{\psi_1}$ with $c_1(j)\ge 0$ and $c_2(j)$
rational numbers.
We have $c_2(0)=c_2(2)=0$ and $c_2(1)=-c_2(3)$.
For a `generic' $g$ all $c_2(j)$ will be zero.
\end{Cor}
{\tt Remark}. Using the result going back to Landau that there are $O(x/\sqrt{\log x})$
integers $n\le x$ having only prime divisors $p$ with $p\equiv 3({\rm mod~}4)$ \cite[pp. 641-669]{Landau2}, it is
easily inferred that the number of non-generic integers $g$ with $|g|\le x$ is
$O(x/\sqrt{\log x})$.
\begin{Cor} {\rm (GRH)}.
Suppose that $h=1$. Then
$\delta_g(3,4;1,4)=\delta_g(3,4;3,4)=1/8$ unless $D$ is divisible by 8 and has
no prime divisor congruent to
$1({\rm mod~}4)$, in which case we have
$$\delta_g(3,4;j,4)={1\over 8}+{\rm sgn}(g)A_{\psi_1}{(-1)^{j+1\over 2}\over
8}\prod_{p|{D\over 8}}{2p\over p^3-p^2-p-1}.$$
\end{Cor}
\begin{Cor} {\rm (GRH)}. Let $j$ be odd. Then
$$\delta_g(3,4;j,4)+\delta_{-g}(3,4;j,4)={1\over 4}.$$
\end{Cor}
An alternative proof of 
the latter corollary is obtained on combining Propostion \ref{prop1} with Theorem \ref{main3}.\\

\noindent {\it Proof of Theorem} \ref{main3explicit}. The proof
is easily deduced from Theorem \ref{main3}.\\
\indent If $h$ is even, then $\Delta_g=0$ and the result follows
on noticing that the density of the set of primes $p\le x$ with $p\equiv 3({\rm mod~}4)$ and $(g/p)=1$
equals $(1+{\rm sgn}(g))/4$.\\
\indent  Next assume that $h$ is odd. Then the set appearing in the
formula for $N_g(3,4;j,4)(x)$ given in Theorem \ref{main3} has density $1/4$.
If $8\nmid D$ it follows from Lemma \ref{intersection} that the summation conditions are never
met  and hence $\Delta_g=0$. So we may assume that $8|D$. By Lemma
\ref{intersection} the $v$'s appearing in the two summations are 
divisible by $D/8$ and thus 
if $D$ has a prime divisor $p$ with $p\equiv 1({\rm mod~}4)$, then
$h_{\psi_1}(D/8)=0$ by Lemma \ref{vauxiliary} and hence $h_{\psi_1}(v)=0$ for these $v$, which
shows that $\Delta_g=0$.\\
\indent It remains to deal with the case where
$8|D$ and $D$ contains no prime divisor $p$ with $p\equiv 1({\rm mod~}4)$.
Using Lemma \ref{degree} and Lemma \ref{intersection} we infer
that
$$\Delta_g={\rm sgn}(g){(-1)^{D+8\over 16}\over 2}\sum_{2\nmid v,~{D\over
8}|v}{h_{\psi_1}(v)(h,v)\over v\varphi(v)}={\rm sgn}(g){(-1)^{D+8\over 16}\over 2}
C_{\chi}(h,2,{D\over 8}).$$
On applying Lemma \ref{eulerproduct} with
$r=2$ and $s=D/8$ (note that since $D$ is a
discriminant $D/8$ must be squarefree), we obtain that
$$\Delta_g={\rm sgn}(g)(-1)^{D+8\over 16}
{P_1P_2P_3\over 2}\prod_{p|D\atop p\equiv 3({\rm
mod~}4)}(-1)=-{\rm sgn}(g){P_1P_2P_3\over 2}.$$
An easy analysis shows that the local factors in the
products $P_1,~P_2$ and
$P_3$ are all
non-negative. \qed

\section{Modulus 3}
\label{modulusdrie}
The case $d=3$ can be dealt with along the
lines of the case $d=4$, hence we supress most details of the proofs.
Our starting point is the following analog of Theorem \ref{triplesum}. 
\begin{Thm}
\label{triplesumtwo}
{\rm (GRH)}. Let $\xi_0,\xi_1$ be the principal, 
respectively non-principal character modulo $3$. For $r\ge 1$
let $\sigma'_{1,r},\sigma'_{-1,r}
\in {\rm Gal}(\mathbb Q(\zeta_{3^{r+1}})/\mathbb Q)$ be the automorphisms
that are uniquely determined by
$\sigma_{1,r}(\zeta_{3^{r+1}})=\zeta_{3^{r+1}}^{1+3^r}$, respectively
$\sigma_{-1,r}(\zeta_{3^{r+1}})=\zeta_{3^{r+1}}^{1+2\cdot 3^r}$.
For $j=-1$ and $j=1$ let 
$$
c'_j(r,tn)=
\cases{1 &if $\sigma'_{j,r}|_{\mathbb Q( \zeta_{3^{r+1}})\cap \mathbb
Q(\zeta_{3^rtn},g^{1/3^rtn})}=$id;\cr
0 &otherwise.}
$$
Let $s\ge 1$. For $j=1$ and $j=2$ we have
$$N_g(1,3^s;j,3)(x)=\delta_g(1,3^s;j,3){x\over \log x}+O\left(\log|g_1g_2|{x\over 
\log^{3/2}x}\right),$$
where
$$\delta_g(1,3^s;j,3)=\sum_{r\ge s}^{\infty}\sum_{t=1\atop 3\nmid t}^{\infty}
\sum_{n=1}^{\infty}{\mu(n)c'_{\xi_1(jt)}(r,tn)
\over [\mathbb Q(\zeta_{3^{r+1}},\zeta_{3^rtn},g^{1/3^rtn}):\mathbb Q]},$$
and the implied constant is absolute.
\end{Thm}
It is very easy to see that $c'_1(r,tn)=c'_{-1}(r,tn)$ for $r\ge 1$ and
thus for $s\ge 1$ we infer that, 
under GRH, $\delta_g(1,3^s;1,3)=\delta_g(1,3^s;2,3)$.
Since (unconditionally)  
$$
\delta_g(1,3^s;0,3)={3^{1-s}\over 2}
-\sum_{r\ge s}\Big({1\over [K_{3^r,3^r}:\mathbb
Q]}-{1\over [K_{3^{r+1},3^r}:\mathbb Q]}\Big),
$$
we then easily deduce, using Lemma \ref{degree}, the following result.
\begin{Thm} {\rm (GRH)}. 
\label{explicietenvolledig}
Let $e_3=\nu_3(h)$ and $s\ge 1$.
If $e_3\le s$, then
$$\cases{\delta_g(1,3^s;0,3)=3^{1-s}/2-3^{2+e_3-2s}/8;\cr
\delta_g(1,3^s;1,3)=3^{2+e_3-2s}/16;\cr
\delta_g(1,3^s;2,3)=3^{2+e_3-2s}/16.}$$
If $e_3>s$, then
$$\cases{\delta_g(1,3^s;0,3)=3^{1-e_3}/8;\cr
\delta_g(1,3^s;1,3)=3^{1-s}/4-3^{1-e_3}/16;\cr
\delta_g(1,3^s;2,3)=3^{1-s}/4-3^{1-e_3}/16.}$$
\end{Thm}
\indent The reason that we cannot take $s=0$ in Theorem \ref{triplesumtwo}
is that $\sigma'_{-1,0}$ does not give rise to an automorphism of
$\mathbb Q(\zeta_3)$. On the other hand $\sigma'_{1,0}$ does and thus we
can define $c'_1(0,tn)$ as in Theorem \ref{triplesumtwo}.\\
\indent Let $j\in \{1,2\}$. 
Since $N_g(0,1;j,3)(x)=N_g(1,3;j,3)(x)+N_g(2,3;j,3)(x)+O(1)$, and
$N_g(1,3;j,3)(x)$ is covered in Theorem \ref{triplesumtwo}, it
remains to deal with $N_g(2,3;j,3)(x)$. Note that
$$N_g(2,3;j,3)(x)=\sum_{t\equiv j({\rm mod~}3)}\#\{p\le x: p\equiv 2({\rm mod~}3),
~r_g(p)=t\}.$$
Reasoning as in Theorem \ref{triplesum}, we then find that, under GRH, we
have
\begin{equation}
\label{bedenkeensiets}
N_g(2,3;j,3)(x)=\delta_g(2,3;j,3){x\over \log x}+
O\left(\log|g_1g_2|{x\over \log^{3/2} x}\right), 
\end{equation}
where
\begin{eqnarray}
\delta_g(2,3;j,3)&=&\sum_{t\equiv j({\rm mod~}3)}\sum_{3\nmid n}
{\mu(n)c'_1(0,tn)\over [\mathbb Q(\zeta_3,\zeta_{tn},g^{1/tn}):\mathbb Q]}\cr
&=&{1\over 2}\sum_{3\nmid v}
{(h_{\xi_0}(v)+\xi_1(j)h_{\xi_1}(v))c'_1(0,v)\over [\mathbb Q(\zeta_3,\zeta_v,g^{1/v}):\mathbb Q]}\cr
&=&{1\over 4}+{\xi_1(j)\over 4}\sum_{3\nmid v\atop \zeta_3\not\in K_{v,v}}
{h_{\xi_1}(v)\over [K_{v,v}:\mathbb Q]},\nonumber
\end{eqnarray}
and the implied constant in (\ref{bedenkeensiets}) is absolute.
In the derivation of the second equality we used Lemma 
\ref{nonmulttomult} and in the derivation of the latter equality we
used the trivial observation that $c_1'(0,v)=1$ iff $\zeta_3\not\in K_{v,v}$.
To sum up we obtained the following theorem.
\begin{Thm}
\label{weetnietmeer}
{\rm (GRH)}. The estimate {\rm (\ref{bedenkeensiets})}
 holds with
$$\delta_g(2,3;j,3)={\xi_0(j)\over 4}+{\xi_1(j)\over 4}
\sum_{3\nmid v\atop \zeta_3\not\in K_{v,v}}{h_{\xi_1}(v)\over 
[K_{v,v}:\mathbb Q]}$$
and an absolute implied constant.
\end{Thm}
\indent (Since the condition $\zeta_3\not\in K_{v,v}$ implies
$3\nmid v$, the latter condition can be dropped, in principle.) 
Theorem \ref{weetnietmeer} is the density version of the following
lemma. 
\begin{Lem}
\label{identity3}
The quantity $N_g(2,3;j,3)(x)$ equals
$${\xi_0(j)\over 2}\#\{p\le x:p\equiv 2({\rm mod~}3)\}+
{\xi_1(j)\over 2}\sum_{3\nmid v}^{\infty}h_{\xi_1}(v)
[\pi_{K_{v,v}}(x)-\pi_{K_{3v,v}}(x)].$$
\end{Lem}
{\it Proof}. Similar to that of Lemma \ref{simpeltoch}. \qed\\

\indent On noting that
$${1\over 2}\sum_{3\nmid v\atop \zeta_3\not\in K_{v,v}}{h_{\xi_1}(v)\over 
[K_{v,v}:\mathbb Q]}=\sum_{3\nmid v}h_{\xi_1}(v)\Big(
{1\over [K_{v,v}:\mathbb Q]}-{1\over [K_{3v,v}:\mathbb Q]}\Big),$$
and invoking Lemma \ref{laatsteh}
and Lemma \ref{eulerproduct}, we obtain the following three
colloraries of Theorem \ref{weetnietmeer}.
\begin{Cor}
\label{driedifference}
{\rm (GRH)}.
Put $\epsilon=1$ if $3\nmid D$ and $\epsilon=-1$ otherwise.
If $g>0$, then 
$$\delta_g(2,3;j,3)={\xi_0(j)\over 4}+{\xi_1(j)\over 4}\Big(C_{\xi_1}(h,3,1)+
\epsilon C_{\xi_1}(h,3,{n_1\over (3,n_1)})\Big).$$
If $g<0$, then 
$\delta_g(2,3;j,3)$ equals $${\xi_0(j)\over 4}+
{\xi_1(j)\over 4}\Big(C_{\xi_1}(h,3,1)
-{C_{\xi_1}(h,3,2)\over 2}+{C_{\xi_1}(h,3,2^{\nu_2(h)+1})\over 2}
+\epsilon C_{\xi_1}(h,3,{n_1\over (n_1,3)})\Big).$$
\end{Cor}
\begin{Cor} {\rm (GRH)}. 
\label{technischhoor}
Recall that $e_p=\nu_p(h)$. Define $\Omega(n)=\sum_{p|n}\nu_p(n)$.
Define
$$P'_1=\prod_{p|D,~p>2\atop p\equiv 2({\rm mod~}3)}
\left({2[p^{e_p}-(-1)^{e_p}]\over p^{e_p-1}(p^2-1)}
+{2p(-1)^{e_p}\over p^{e_p}(p^2+1)(p-1)}\right){\rm ~and}$$
$$P'_2=\prod_{p\nmid 2D\atop p\equiv 2({\rm mod~}3)}
\left(1-{2[p^{e_p}-(-1)^{e_p}]\over p^{e_p-1}(p^2-1)}
+{2p(-1)^{e_p+1}\over p^{e_p}(p^2+1)(p-1)}\right).$$
Put $\epsilon_1=0$ if $D$ has a prime divisor
$q$ that satisfies $q\equiv 1({\rm mod~}3)$ and $\epsilon_1=1$
otherwise.  
If $g>0$, then
$$\delta_g(2,3;j,3)={\xi_0(j)\over 4}+{\xi_1(j)\over 4}C_{\xi_1}(h,3,1)+
\epsilon_1{\xi_1(j)\over 5}(-1)^{\Omega(n_1)}2^{e_2+2-2\nu_2(n_1)}P'_1P'_2.$$
If $g<0$, then
$$\delta_g(2,3;j,3)={\xi_0(j)\over 4}+{\xi_1(j)\over 4}
\Big(1-{[2^{e_2}-(-1)^{e_2}]\over 
3\cdot 2^{e_2-1}}-{2^{2-e_2}(-1)^{e_2}\over 5}\Big)C_{\xi_1}(h,6,1)+$$
$$\epsilon_1{\xi_1(j)\over 5}(-1)^{\Omega(n_1)}2^{e_2+2-2\nu_2(n_1)}P'_1P'_2.$$
\end{Cor}
\begin{Cor} {\rm (GRH)}.
Suppose $h=1$. We have
$$\delta_g(2,3;j,3)={\xi_0(j)\over 4}+{\xi_1(j)\over 4}A_{\xi_1}
\Big(1+\epsilon_1(-1)^{\Omega(n_1)}2^{4-2\nu_2(n_1)}
\prod_{p|D\atop p>3}{2p\over p^3-p^2-p-1}\Big).$$
\end{Cor}
{\tt Remark}. Note that
$$C_{\xi_1}(h,3,1)=\prod_{p\equiv 2({\rm mod~}3)}
\left(1-{2[p^{e_p}-(-1)^{e_p}]\over p^{e_p-1}(p^2-1)}
+{2p(-1)^{e_p+1}\over p^{e_p}(p^2+1)(p-1)}\right).$$

\noindent A somewhat tedious analysis of Corollary \ref{technischhoor}
together with Theorem \ref{explicietenvolledig} 
yields the
following size comparison of $\delta_g(2,3;1,3)$ with $\delta_g(2,3;2,3)$ 
and of $\delta_g(1,3)$ with $\delta_g(2,3)$.
\begin{Prop} {\rm (GRH)}. 
If $g>0$ and $h$ is even, then $\delta_g(2,3;1,3)\le \delta_g(2,3;2,3)$,
otherwise $\delta_g(2,3;1,3)\ge \delta_g(2,3;2,3)$. We have
$\delta_g(2,3;1,3)=\delta_g(2,3;2,3)$ iff 
$\mathbb Q(\sqrt{g_0})=\mathbb Q(\sqrt{3})$ and
$\nu_2(h)\in \{0,2\}$.
The same result holds
with $\delta_g(2,3;j,3)$ replaced by $\delta_g(j,3)$.
\end{Prop}

\section{On the generic behaviour of $\delta_g(a,d)$, $d=3,4$}
If $g$ is not a square or -1, then an old heuristic model predicts that
the number of primes $p\le x$ such that $g$ is a primitive root mod $p$
should be asymptotically equal to $\sum_{p\le x}\varphi(p-1)/(p-1)$, where
$\varphi(p-1)/(p-1)$ is the density of primitive roots in $\mathbb F_p^*$.
It is easily proved, 
see e.g. \cite{Heuristics}, that on average $\varphi(p-1)/(p-1)$ is equal to the
Artin constant $A$, that is
$$\lim_{x\rightarrow \infty}{1\over \pi(x)}\sum_{p\le x}{\varphi(p-1)\over p-1}
=A=0.37395\cdots.$$ 
From the work of Hooley 
\cite{Hooley} it can be deduced that under GRH for a positive 
proportion of all $g$ the above heuristic is false.\\
\indent Let $\delta(p;a,d)=\sum_{r|p-1,~r\equiv a({\rm mod~}d)}\varphi(r)/(p-1)$,
then $\delta(p;a,d)$ is the density of elements in the multiplicative group of
the finite field $\mathbb F_p$ with order congruent to $a({\rm mod~}d)$.
A (naive) heuristic prediction for $N_g(a,d)(x)$ is then provided by
$\sum_{p\le x}\delta(p;a,d)$. It can be shown that
$\lim_{x\rightarrow \infty}\sum_{p\le x}\delta(p;a,d)/\pi(x)=\delta(a,d)$
exists \cite{Moreeaverage}. For $d=3,4$ some computation \cite{Moreeaverage} shows that
$$\delta(a,3)=\cases{{3\over 8} &if $a\equiv 0({\rm mod~}3)$;\cr
{5\over 16}+{A_{\xi_1}\over 4} & if $a\equiv 1({\rm mod~}3)$;\cr
{5\over 16}-{A_{\xi_1}\over 4} &if $a\equiv 2({\rm mod~}3)$}$$
and
$$\delta(a,4)=\cases{{1\over 3} &if $a$ is even;\cr
{1\over 6} &if $a$ is odd.}$$
On comparing this computation with our conditional results for
$\delta_g(a,3)$ and $\delta_g(a,4)$ we obtain the following result.
\begin{Prop} {\rm (GRH)}. Let $d=3,4$ be fixed. There are
at most $O(x/\sqrt{\log x})$ integers $|g|\le x$ for which
$\delta_g(a,d)\ne \delta(a,d)$ for some integer $a$.
In particular, for almost all integers $|g|\le x$ we have
$\delta_g(a,d)=\delta(a,d)$ for every integer $a$.
\end{Prop}
This proposition shows
that for fixed $d=3,4$ it makes sense to call an integer $g$ 
{\it generic}
if $\delta_g(a,d)=\delta(a,d)$ for every integer $a$.\\
\indent In a similar vein we have:
\begin{Prop} {\rm (GRH)}. Let $d=3,4$ be fixed.
If $|D(g)|$ tends to infinity as $g$ ranges over a set of
rationals $g$ for which $h=1$, then $\delta_g(a,d)$ tends to $\delta(a,d)$.
\end{Prop}
The latter two  results seem to hold for other values of $d$ as well (cf. Table 2 of
\cite{Moreeaverage}), with $O(x\/\sqrt{\log x})$ replaced by $o(x)$. I might return to
this in a sequel.\\

\noindent Acknowledgement.  The work presented here was carried out at the 
Korteweg-de Vries Institute (Amsterdam) in the NWO Pioneer-project of Prof. E. Opdam.
I am grateful to Prof. Opdam for the opportunity to work in his group.\\
\indent I thank K. Chinen, L. Murata, H. Roskam and P. Tegelaar for
pointing out some inaccuracies in earlier versions and K. Chinen
for showing me some of his (extensive) numerical work. The data given
in Tables 1 and 2 were calculated using a $C^{++}$ program kindly written by Yves Gallot. 
Computations took on average (per case) 25 minutes on a 1.50 GHz Pentium IV processor.

\vfil\eject
\section{Tables}
We illustrate our results by some examples (assuming GRH).\\
\hfil\break
\centerline{{\bf Table 1:} The case $d=3$}
\begin{center}
\begin{tabular}{|c|c|c|c|c|c|}\hline
$g$&$g_0$&$h$&$\delta_g(1,3)-\delta_g(2,3)$&numerical&experimental\\
\hline\hline
$-14^{4}$&14&4&$3A_{\xi_1}/4$&$+0.13048284\dots$&$+0.13045317$\\ \hline
-196&14&2&$A_{\xi_1}$&$+0.17397712\dots$&$+0.17399131$\\ \hline
$-3^8$&3&8&$15A_{\xi_1}/16$&$+0.16310355\dots$&$+0.16310903$\\ \hline
-3&3&1&$5A_{\xi_1}/2$&$+0.43494280\dots$&$+0.43499017$\\ \hline
-2&2&1&$3A_{\xi_1}/8$&$+0.06524142\dots$&$+0.06525031$\\ \hline
3&3&1&0&0&$+0.00001393$\\ \hline
9&3&2&$-5A_{\xi_1}/2$&$-0.43494280\dots$&$-0.43502303$\\ \hline
81&3&4&0&0&$-0.00001895$\\ \hline
6561&3&8&$-5A_{\xi_1}/4$&$-0.21747140\dots$&$-0.21748481$\\ \hline
2&2&1&$3A_{\xi_1}/8$&$+0.06524142\dots$&$+0.06515583$\\ \hline
4&2&2&$-7A_{\xi_1}/4$&$-0.30445996\dots$&$-0.30442279$\\ \hline
5&5&1&$67A_{\xi_1}/94$&$+0.12400497\dots$&$+0.12397327$\\ \hline
25&5&2&$-151A_{\xi_1}/94$&$-0.27947388\dots$&$-0.27952119$\\ \hline
49&7&2&$-3A_{\xi_1}/2$&$-0.26096568\dots$&$-0.26097396$\\ \hline
2401&7&4&$-A_{\xi_1}/2$&$-0.08698856\dots$&$-0.08697494$\\ \hline
\end{tabular}
\end{center}
\medskip
\centerline{{\bf Table 2:} The case $d=4$}
\begin{center}
\begin{tabular}{|c|c|c|c|c|c|}\hline
$g$&$g_0$&$h$&$\delta_g(1,4)-\delta_g(3,4)$&numerical&experimental\\
\hline\hline
-216&6&3&$9A_{\psi_1}/28$&$+0.20688771\dots$&$+0.20686925$\\ \hline
-9&3&2&0&0&$+0.00000068$\\ \hline
-81&3&4&0&0&$-0.00000232$\\ \hline
2&2&1&$-A_{\psi_1}/4$&$-0.16091266\dots$&$-0.16088852$\\ \hline
4&2&2&0&0&$+0.00001122$\\ \hline
8&2&3&$-A_{\psi_1}/28$&$-0.02298752\dots$&$-0.02301736$\\ \hline
512&2&9&$-3A_{\psi_1}/28$&$-0.06896257\dots$&$-0.06897632$\\ \hline
216&6&3&$-9A_{\psi_1}/28$&$-0.20688771\dots$&$-0.20687020$\\ \hline
2048&2&11&$-489A_{\psi_1}/2396$&$-0.13136276\dots$&$-0.13134226$\\ \hline
$6^9$&6&9&$-A_{\psi_1}/4$&$-0.16091266\dots$&$-0.16088478$\\ \hline
$6^{27}$&6&27&$-23A_{\psi_1}/84$&$-0.17623768\dots$&$-0.17620628$\\ \hline
\end{tabular}
\end{center}
\hfil\break
\hfil\break
The number in the column `experimental' arose on taking the density difference
over the first $10^8$ primes, not letting the primes $p$ for which the order
of $g$ mod $p$ is not defined contribute to either $\delta_g(1,3)$ or
$\delta_g(2,3)$ (in Table 1), or $\delta_g(1,4)$ and $\delta_g(3,4)$ in
Table 2. Thus, for example, in the column headed `experimental' in Table 1
the numbers 
$${N_g(1,3)(p_{10^8})-N_g(2,3)(p_{10^8})\over 10^8}$$
are recorded (recall that $p_{10^8}=2038074743$).
The last decimals in the columns headed `numerical' and `experimental' are not rounded.\\

\end{document}